\documentclass[12pt]{amsart}

\usepackage{amscd,amssymb,amsthm,verbatim}

\newcommand{\C}{{\mathbb C}}

\newcommand{\const}{\operatorname{const.}}

\newcommand{\diam}{\operatorname{diam}}

\newcommand{\dvol}{\operatorname{dvol}}

\newcommand{\HH}{\operatorname{H}}

\newcommand{\Id}{\operatorname{Id}}

\newcommand{\Image}{\operatorname{Im}}

\newcommand{\pt}{\operatorname{pt}}

\newcommand{\R}{{\mathbb R}}

\newcommand{\re}{\operatorname{Re}}

\newcommand{\Ric}{\operatorname{Ric}}
\newcommand{\Rm}{\operatorname{Rm}}

\newcommand{\sing}{\operatorname{sing}}

\newcommand{\Tr}{\operatorname{Tr}}

\newcommand{\vol}{\operatorname{vol}}
\newcommand{\Vol}{\operatorname{Vol}}
\newcommand{\Z}{{\mathbb Z}}

\newcommand{\eps}{\epsilon}

\numberwithin{equation}{section}
\theoremstyle{plain}

\newtheorem{lemma}[equation]{Lemma}
\newtheorem{theorem}[equation]{Theorem}

\newtheorem{corollary}[equation]{Corollary}

\theoremstyle{definition}
\newtheorem{definition}[equation]{Definition}

\theoremstyle{definition}
\newtheorem{example}[equation]{Example}

\theoremstyle{definition}
\newtheorem{remark}[equation]{Remark}
\def\<{\langle}
\def\>{\rangle}
\def\({\left(}
\def\){\right)}
\def\p{\partial}

\begin{document}

\title{Ricci flow on quasiprojective manifolds}
\author{John Lott and Zhou Zhang}
\address{Department of Mathematics\\
University of California - Berkeley\\
Berkeley, CA  94720-3840\\ 
USA} \email{lott@math.berkeley.edu}
\address{Department of Mathematics\\
University of Michigan\\
Ann Arbor, MI  48109-1109\\ 
USA}
\email{zhangou@umich.edu}

\thanks{The research of the first author was partially supported
by NSF grant DMS-0903076}
\thanks{The research of the second author was partially supported
by NSF grant DMS-0904760}
\date{November 18, 2010}

\begin{abstract}
We consider the K\"ahler-Ricci flow on complete 
finite-volume metrics that
live on the complement of a divisor
in a compact K\"ahler manifold $\overline{X}$.  Assuming certain
spatial asymptotics on the initial metric, we compute the
singularity time in terms of cohomological data on
$\overline{X}$. We also give a sufficient condition for the
singularity, if there is one, to be type-II.
\end{abstract}

\maketitle

\tableofcontents

\section{Introduction}

In this paper we study the Ricci flow on certain finite-volume complete
K\"ahler metrics that live on complements of divisors
in compact K\"ahler manifolds. Our motivation, which we now describe,
comes from the general goal of
understanding singularities in Ricci flow.

It is known, since Hamilton's first Ricci flow paper
\cite{Hamilton}, that 
singularities in a Ricci flow on a manifold $M$
arise from curvature blowup.
The nature of the blowup is important in the analysis
of the singularity. We let $\Rm$ denote sectional curvatures.
If $T_{\sing}$ is a first singularity time then the singularity is
said to be {\em type-I} if there is a constant $C < \infty$ so that
$|\Rm|(m,t) \le \frac{C}{T_{\sing}-t}$ for all $m \in M$ and
$t < T_{\sing}$. Otherwise,
the singularity is said to be {\em type-II}.

In Ricci flow, the natural scaling is that $time \sim
distance^2$.  Since $curvature \sim distance^{-2}$,
a naive dimensional analysis would suggest that all
singularities are type-I.  However, this is not the case.
The first type-II singularity was
found on a noncompact surface which is diffeomorphic to $\R^2$, but
whose initial metric $g(0)$ describes a hyperbolic
cusp capped off by a ball.  The singular time is
$T_{\sing} = \frac{1}{4\pi} \Vol(\R^2, g(0))$.
At any time $t < T_{\sing}$, the volume is $\Vol(\R^2, g(0)) - 4 \pi t$.
Hence as $t \rightarrow T_{\sing}$, there is no volume left.
The geometric behavior as the time $t$ approaches $T_{\sing}$ is as follows.
For $t < T_{\sing}$,
one can divide the surface into an inner region $I_t$ and
an outer region $O_t$. 
As one goes out the end, the metric on $O_t$ has asymptotically constant 
negative curvature $k(t)$, with $k(t)$ remaining bounded as
$t \rightarrow T_{\sing}$. However, as $t \rightarrow T_{\sing}$, the
outer region disappears and the inner region
$I_t$ dominates. The curvature on $I_t$ goes to infinity
pointwise as
$t \rightarrow T_{\sing}$ and its geometry approaches a ray, in the
pointed Gromov-Hausdorff sense.
After a parabolic rescaling to normalize the curvature, $I_t$ approaches
a special Ricci flow solution, the cigar soliton,
as one approaches the singularity time.
For these results and more, we refer to papers by
Daskalopoulos-del Pino-Hamilton-Sesum 
\cite{Daskalopoulos-del PinoI,Daskalopoulos-del PinoII,Daskalopoulos-Hamilton,Daskalopoulos-Sesum} and references therein.

The goal of this paper is to extend some of these two-dimensional results
to higher dimensions.  A starting point in the two-dimensional analysis
is the use of isothermal coordinates on $\R^2$, 
in order to write the Ricci flow
equation as a scalar equation for the conformal factor.  This method
clearly does not work in higher dimensions, so we must take another approach.
Our approach is based on the 
observation that $\R^2$, with a finite-volume asymptotically
hyperbolic metric, can be considered as $S^2 - \pt$, with a metric on
$S^2$ which, in local coordinates near $\pt$, approaches the
Poincar\'e metric $\frac{4 dz d\overline{z}}{|z|^2 \log^2
(|z|^{-2})}$. This is an example of a
quasiprojective manifold, meaning the complement $X = \overline{X} - D$
of an effective divisor $D$ with simple normal crossings 
in a compact K\"ahler manifold 
$\overline{X}$. Another simple example of a quasiprojective manifold
comes from taking the
product $X = (S^2 - \pt) \times (S^2 - \pt)$ of the previous manifold
with itself. Then $X = \overline{X} - D$, where
$\overline{X} = S^2 \times S^2$ and $D = (S^2 \times \pt) \cup (\pt \times S^2)$.

In what follows, we will speak equivalently of a K\"ahler metric or a
K\"ahler form.
Let $\omega_X(0)$ be a complete K\"ahler metric with bounded sectional
curvature on a complex manifold $X$. It is known that
there is some $\eps > 0$ so that there is a 
Ricci flow solution on the time interval $[0,\epsilon]$ with
initial metric $\omega_X(0)$,
complete time slices and uniformly bounded sectional curvature
\cite{shi}.
It is easy to see that the time-$t$ metric $g(t)$ is K\"ahler with
respect to the initial (and fixed) complex structure, so it makes sense
to talk about the ensuing K\"ahler-Ricci flow
\cite{cao}.

By definition, the singularity time $T_{\sing}$ 
is the supremum of
the numbers $T > 0$ with the property that
there is a Ricci flow solution
$\omega_X(t)$ with the given value at $t=0$, defined for
$t \in [0, T]$, having complete time slices and
uniformly bounded sectional curvature on
the time interval $[0,T]$. Note that $T_{\sing}$ could be infinity,
which corresponds to not having a singularity.

In order to state the main result, we introduce some terminology.
Given a compact K\"ahler manifold $\overline{X}$ of complex dimension $n$, we write
$[K_{\overline{X}}] \in \Image \left(\HH^2(\overline{X}; \Z) \rightarrow
\HH^2(\overline{X}; \R) \right) \cap \HH^{(1,1)}(\overline{X}; \R)$
for the first Chern class of the canonical line bundle 
$K_{\overline{X}} = \Lambda^{n,0}\overline{X}$.
Note that $[K_{\overline{X}}]$ is the negative of the first Chern class of the holomorphic
tangent bundle, so $[K_{\overline{X}}] = - c_1(\overline{X})$.

For us,
a divisor $D = \sum_i D_i$ in $\overline{X}$ is a formal sum of 
closed complex submanifolds
of $\overline{X}$ with complex codimension one.
There is a corresponding class 
$[D] \in \Image \left(\HH^2(\overline{X}; \Z) \rightarrow
\HH^2(\overline{X}; \R) \right) \cap \HH^{(1,1)}(\overline{X}; \R)$,
whose Poincar\'e dual $*[D] \in \HH_{2n-2}(\overline{X}; \R)$
is the sum of the pushforwards of the fundamental classes of the $D_i$'s.
Hereafter we will assume that $D$ has normal crossings.

Recall that a class $c \in \HH^{(1,1)}(\overline{X}; \R)$ is a {\em K\"ahler class}
if there is a closed positive form $\omega \in \Omega^{(1,1)}(\overline{X})$ whose
de Rham cohomology class is $c$. In such a case, we write $c > 0$.

The main theorem of the paper concerns the 
K\"ahler-Ricci flow solution on $X = \overline{X} - D$ whose initial
metric is a finite-volume K\"ahler metric $\omega_X(0)$ with ``superstandard''
spatial asymptotics.  This notion, which will be made precise in Definition
\ref{superstandard}, roughly means that the metric at infinity can be
decomposed into
families of products of hyperbolic cusp metrics.  (An example of superstandard
spatial asymptotics is the product metric on $X = (S^2 - \pt) \times (S^2 - \pt)$
from before.) One motivation for considering such asymptotics is that
they arise for the finite-volume
K\"ahler-Einstein metric on $X$ that exists when
$[K_{\overline{X}}+D] > 0$ \cite{kobayashi,schumacher,tian-yau,tsu,wu}.

If $\omega_X(0)$ has superstandard spatial asymptotics then in terms of
the inclusion $X \subset \overline{X}$, we can extend $\omega_X(0)$ by zero to
get a closed $(1,1)$-current on $\overline{X}$; see
Theorem \ref{extension}. There is a corresponding
cohomology class $[\omega_X(0)] \in \HH^{(1,1)}(\overline{X}; \R)$.

The goal now is to express properties of the K\"ahler-Ricci flow on $X$ in
terms of cohomological data. We remark that 
it may not be immediately clear which 
cohomology group is the relevant one. For example, one may think that
it should be some sort of cohomology of $X$. 
However, it turns out that what's relevant is the cohomology of
the compactification $\overline{X}$.
(As a precedent, the cohomology of the compactification is also
key to the previously-mentioned
work on finite-volume K\"ahler-Einstein manifolds.)
We show
that we can effectively compute $T_{\sing}$ from cohomological data
on $\overline{X}$. We also give a sufficient condition to ensure a
type-II singularity.

\begin{theorem} \label{maintheorem}
Suppose that $\omega_X(0)$ is a K\"ahler metric on $X$ with superstandard
spatial asymptotics. \\
1. The singularity time $T_{\sing}$ of the ensuing 
(unnormalized) K\"ahler-Ricci flow
equals the supremum of
the numbers $T>0$ so that $[\omega_X(0)] + 2\pi T [K_{\overline{X}}+D] \in 
\HH^{(1,1)}(\overline{X}; \R)$ is a K\"ahler class
on $\overline{X}$. \\
2. If $D \neq \emptyset$, $T_{\sing} < \infty$ and
$[\omega_X(0)] + 2\pi T_{\sing} [K_{\overline{X}}+D]$ vanishes in 
$\HH^{(1,1)}(\overline{X}; \R)$ then there is a type-II singularity at time
$T_{\sing}$.
\end{theorem}

When $X$ has one complex dimension, Theorem
\ref{maintheorem} recovers some of the surface results mentioned before;
see Example \ref{surfaceexample}.

In the course of proving Theorem \ref{maintheorem} we obtain some
results about K\"ahler-Ricci flow that are valid for a wider
class of initial metrics. We now describe some of these results, in
order of decreasing generality. 

In Theorem \ref{th:existence1} we characterize
the singularity time for a
normalized K\"ahler-Ricci flow on any complex manifold $X$, 
whose initial metric $\omega_X(0)$ is complete with bounded
curvature. For $T \ge 0$, put $\omega_T = 
-{\rm Ric}(\omega_X(0))+e^{-T}\left(\omega_X(0)
+{\rm Ric}(\omega_X(0))\right)$.
Theorem \ref{th:existence1} says that
the singularity time $T_{\sing}$ equals the supremum of the
numbers $T > 0$ with the property that there is some
$F_T \in C^\infty(X)$ so that  
\begin{itemize}
\item $\omega_T + \sqrt{-1} \partial
\overline{\partial}F_T$ is a K\"ahler metric on $X$ which is
biLipschitz to $\omega_X(0)$, and
\item For each $k \ge 0$, the $k$-th
covariant derivatives of $F_T$ (with respect to the initial
metric $\omega_X(0)$) are uniformly bounded on $X$.
\end{itemize}
Theorem \ref{th:existence1} 
is an extension of \cite[Proposition 1.1]{t-znote} by
Tian and Zhang, which 
dealt with the case when $X$ is compact. The interest of Theorem
\ref{th:existence1} is that the issue of computing $T_{\sing}$ is
reduced to a flow-independent question on $X$.

The next main result,
Theorem \ref{th:non-deg-convergence1}, concerns long-time convergence.
Under the assumption that the initial metric $\omega_X(0)$ satisfies
$-{\rm Ric}(\omega_X(0))+\sqrt{-1}\partial\bar\partial 
f>\epsilon\omega_X(0)$ for some $\epsilon > 0$ and some 
smooth function $f$ with bounded 
covariant derivatives, we show that the normalized K\"ahler-Ricci flow
(\ref{eq:krf}) exists forever and that its time slices
converge smoothly to a complete K\"ahler-Einstein metric on $X$
with Einstein constant $-1$.  This is an extension of 
\cite[Theorem 1.1]{chau} by Chau.

The next goal is to characterize the singularity time in
cohomological terms.
To do so, we specialize to initial metrics on 
a quasiprojective manifold
$X = \overline{X} - D$ that satisfy
``standard'' spatial asymptotics. 
In Theorem \ref{0thorder} we show that this property
is shared by the time slices of the ensuing 
normalized K\"ahler-Ricci flow. 
We can then extend $\omega_X(t)$ by zero to define a
closed $(1,1)$ current on $\overline{X}$ and a
corresponding cohomology class $[\omega_X(t)] \in
\HH^{(1,1)}(\overline{X}; \R)$.
We prove that $[\omega_X(t)]$ equals
$e^{-t} [\omega_X(0)] \: + \: 2\pi (1 - e^{-t}) \: 
[K_{\overline{X}} + D]$.
In Theorem \ref{extension}
we show that if the K\"ahler-Ricci flow on $X$,
with initial metric $\omega_X(0)$, extends to time $T$ then
$[\omega_X(T)]$ is a K\"ahler class on $\overline{X}$. The proof
uses a characterization of K\"ahler classes that is due to
Demailly-Paun \cite{demailly-paun}.

In Theorem \ref{superextends} we further specialize to initial metrics on 
$X = \overline{X} - D$ that satisfy
``superstandard'' spatial asymptotics.  We show that this property
is again shared by the time slices of the ensuing K\"ahler-Ricci flow.
In Theorem \ref{end1}
we show that if $[\omega_X(t)]$ happens to be a K\"ahler class on
$\overline{X}$ then $\omega_X(t)$ can be written as
$\omega_t + \sqrt{-1} \partial
\overline{\partial}F_t$ for an appropriate $F_t \in C^\infty(X) \cap
L^\infty(X)$. Along with Theorem \ref{th:existence1}, this proves the
first part of Theorem \ref{maintheorem}.

The proof of the second part of Theorem \ref{maintheorem} is by
contradiction.  Suppose that the singularity is type-I.
By a result of Naber \cite{Naber} (which is based on Perelman's work
\cite{Perelman}),
there is a spacetime sequence $(x_i,t_i)$ with 
$t_i \rightarrow T_{\sing}$ so that after rescaling by
$\frac{1}{T_{\sing} - t_i}$, the corresponding pointed
Ricci flow solutions converge to a $\kappa$-noncollapsed gradient
shrinking soliton $Y$ with uniformly bounded curvature. 
Since $D \neq \emptyset$, the manifold $Y$ is noncompact. By our
assumption on the limit of the K\"ahler class, $Y$ has finite volume.
This leads to a contradiction.  Therefore, the singularity 
must be type-II.

We mention some open problems. The first problem is to understand
what kind of rescaling limits can arise from type-II singularities
as above.
A general construction of Hamilton gives a rescaling limit
which is an eternal solution, i.e. which exists for $t \in \R$
\cite[Proposition 8.17]{Chowetal}. The question is whether it must be a gradient steady
soliton, as is the case in one complex dimension, where one gets
the cigar soliton. Another question
is which gradient steady solitons can occur as rescaling limits.

A second problem is to work out a spatial asymptotic expansion for
the metric $\omega_X(t)$, assuming some precise spatial asymptotics for
$\omega_X(0)$. The analogous question for a K\"ahler-Einstein metric
on $X = \overline{X} - D$, which exists when $[K_{\overline{X}} +D] >0$, was
addressed in \cite{schumacher,wu}.

We thank Lei Ni for a helpful comment. We thank the referee for
a careful reading and helpful suggestions.

\section{Conventions}

Given a K\"ahler manifold $X$ of complex dimension $n$,
the K\"ahler form is a real $(1,1)$-form $\omega$ 
which can be expressed in holomorphic normal coordinates at a 
point $p$ by $\omega(p) = \frac{\sqrt{-1}}{2} 
\sum_{i=1}^n dz^i \wedge d\overline{z}^i$. The 
K\"ahler form of the Poincar\'e metric is given on 
the upper half plane $H = \{ w \in \C : \text{Im}(w) > 0\}$ 
by $\frac{\sqrt{-1}}{2} \frac{dw \wedge d\overline
{w}}{(\text{Im}(w))^2}$. This is the pullback of the K\"ahler 
form
\begin{equation}
2 \sqrt{-1} \frac{dz \wedge d\overline{z}}{|z|^2 \log^2
(|z|^{-2})} = - \sqrt{-1} \partial \bar{\partial} \log 
\left( |z|^2 \log^2(|z|^{-2})\right)
\end{equation} 
on 
$\Delta^* = \Delta - \{0\}$, under the map $z = e^{\sqrt{-1}w}$.
Here $\Delta$ denotes the unit ball in $\C$.

Let $L$ be a holomorphic line bundle over $X$ with 
Hermitian metric $h_L$. If $\sigma$ is a section of $L$ then
we write $|\sigma|^2_L =
h_L(\sigma, \sigma)$. 
There is a unique connection $\nabla^L$ which is compatible
with both the Hermitian metric $h_L$ and the holomorphic
structure on $L$.
Let $F(h_L) \in \Omega^2(M)$ 
be the curvature form of $\nabla^L$.
The de Rham cohomology class of $\frac{\sqrt{-1}}{2\pi} F(h_L)$ equals
$c_1(L) \in \Image(\HH^2(X; \Z) \rightarrow \HH^2(X; \R))$. 
If $\sigma$ is a local holomorphic section 
of $L$ then $F(h_L) \: = \: - \: \partial \overline
{\partial} \log |\sigma|^2_{h_L}$.
If $K_X = \Lambda^{n,0}X$ is the canonical bundle 
of $X$ then we will write $[K_X]$ for $c_1(K_X) \: = \: - \: c_1(X)$.

The Ricci form is 
\begin{equation}
\Ric = - \sqrt{-1} F \left( h_{K_X} \right) = \sqrt{-1} \partial \overline
{\partial}\log |\sigma|_{K_X}^2 = - \sqrt{-1} \partial \overline
{\partial}\log \det(g_{i \overline{j}}), 
\end{equation}
where $\sigma$ is locally $dz^1 \wedge \ldots \wedge dz^n$.
Then $[\Ric] = - 2\pi c_1(K_X) = 2\pi c_1(X) \in \HH^2(X; \R)$. 
For the Poincar\'e metric on $\C^*$, $|\sigma|_{K_X}^2 = 
|z|^2 \log^2(|z|^{-2})$, so $\Ric (\omega)= - \omega$.

\section{The potential flow}

We consider Ricci flow on a connected complex manifold 
$X$ of complex dimension $n$, which may be non-compact. 
Suppose that $\omega_0$ is a smooth complete K\"ahler 
metric on $X$ with bounded curvature. The unnormalized
K\"ahler-Ricci flow equation is
\begin{equation} \label{unnorm}
\frac{\partial\widetilde\omega_t}{\partial t}=-{\rm Ric}
(\widetilde\omega_t),~~~~ \widetilde\omega_0
=\omega_0,
\end{equation}
while for us the normalized K\"ahler-Ricci flow equation is
\begin{equation}
\label{eq:krf}
\frac{\partial\widetilde\omega_t}{\partial t}=-{\rm Ric}
(\widetilde\omega_t)-\widetilde\omega_t,~~~~ \widetilde\omega_0
=\omega_0.
\end{equation}
(Note that the right-hand side of (\ref{unnorm}) differs by a
factor of two from the usual Ricci flow equation
$\frac{dg}{dt} = - 2 \Ric$.) One can pass between
solutions of (\ref{unnorm}) and (\ref{eq:krf}) by rescaling the metric
and reparametrizing time, so there is no essential difference
between the two equations.  Theorem \ref{maintheorem}
is stated for the unnormalized equation (\ref{unnorm}) but in the
rest of the paper, we use the normalized equation 
(\ref{eq:krf}). The reason is that the Poincar\'e metric is
a static solution of (\ref{eq:krf}); this fact will be convenient
in some constructions.

There is some 
$T > 0$ so that there is a solution of 
(\ref{eq:krf})
on the time interval $[0,T]$, having complete time slices 
and uniformly bounded curvature on $[0, T]$ \cite{shi}. Furthermore, 
such a solution is unique on $[0,T]$ \cite{chen-zhu}. 

As is standard in K\"ahler-Ricci flow, we reduce (\ref{eq:krf}) to
a scalar equation. To do so, note that if we have a solution of
(\ref{eq:krf}) then after passing to de Rham cohomology, we get
an ordinary differential equation in $\HH^2(X; \R)$ :
\begin{align} \label{ode}
\frac{d}{dt} [\widetilde\omega_t]&=- [{\rm Ric}
(\widetilde\omega_t)]-[\widetilde\omega_t] \\
&=- [{\rm Ric}
(\omega_0)]-[\widetilde\omega_t]. \notag
\end{align}
The solution to (\ref{ode}) is 
$[\widetilde{\omega}_t]=-[{\rm Ric}(\omega_0)]+e^{-t}\left([\omega_0]
+[{\rm Ric}(\omega_0)]\right)$. This suggests putting
\begin{equation}
\label{omegat}
\omega_t=-{\rm Ric}(\omega_0)+e^{-t}\left(\omega_0
+{\rm Ric}(\omega_0)\right)
\end{equation} 
and making an ansatz for a solution of (\ref{eq:krf}), of the form
$\omega_t + \sqrt{-1} \partial \overline{\partial} u$ for some
scalar function $u$.

Consider the equation
\begin{equation}
\label{eq:skrf2}
\frac{\partial u}{\partial t}={\rm log}\frac{\left(\omega_t+
\sqrt{-1}\partial\bar{\partial} u\right)^n}{\omega^n_0}-u, 
~~~~ u(0, \cdot)=0.
\end{equation}
It is implicit that we only consider solutions $u$ of (\ref{eq:skrf2}) 
on time intervals so that $\omega_t+\sqrt{-1}\partial\bar
{\partial} u > 0$. Note that a solution of (\ref{eq:skrf2}) 
automatically has $\frac{\partial u}{\partial t} (0, \cdot)=
0$.

\begin{lemma}

Suppose that there is a solution to (\ref{eq:krf}) on a time 
interval $[0,T]$, with complete time slices and uniformly 
bounded curvature.  Then there is a smooth solution for 
$u$ in (\ref{eq:skrf2}) on the time interval $[0,T]$ so that \\
1. For each $t\in [0, T]$, $\omega_t+\sqrt{-1}\partial\bar
{\partial} u$ is a K\"ahler metric which is biLipschitz 
equivalent to $\omega_0$. \\
2. For each $k$, the $k$-th covariant derivatives of $u$ 
(with respect to the initial metric $\omega_0$) are 
uniformly bounded. 

Also, $\widetilde{\omega}_t=\omega_t+\sqrt{-1}\partial\bar
{\partial} u$.

Conversely, suppose that there is a smooth solution to 
(\ref{eq:skrf2}) on a time interval $[0,T]$ so that \\
1. For each $t \in [0, T]$, $\omega_t+\sqrt{-1}\partial
\bar{\partial} u$ is a K\"ahler metric which is biLipschitz 
equivalent to $\omega_0$. \\
2. For each $k$, the $k$-th covariant derivatives of $u$ 
(with respect to the initial metric $\omega_0$) are 
uniformly bounded.

Then $\widetilde{\omega}_t=\omega_t+\sqrt{-1}\partial\bar
{\partial} u$ is a solution to (\ref{eq:krf}) on $[0,T]$, with 
complete time slices and uniformly bounded curvature. 

\end{lemma}

\begin{proof}
Suppose that we have a solution to (\ref{eq:krf}) on a time 
interval $[0,T]$, with complete time slices and uniformly 
bounded curvature. Put 
\begin{equation} \label{ueqn1}
u(t) \: = \: \int_0^t e^{s-t} \: {\rm log}\frac{\widetilde{\omega}_s^n}
{\omega^n_0} \: ds,
\end{equation}
so that
\begin{equation} \label{ueqn2}
\frac{\partial u}{\partial t}={\rm log}\frac{\widetilde{\omega}_t^n}
{\omega^n_0}-u.
\end{equation}
Then for each $k$, the $k$-th covariant derivatives of $u$ 
(with respect to the initial metric $\omega_0$) are uniformly 
bounded. Also,
\begin{equation}
\frac{\partial}{\partial t} \left( \widetilde{\omega}_t - \omega_t -  
\sqrt{-1} \partial \bar{\partial} u \right)
= - \left( \widetilde{\omega}_t - \omega_t -  \sqrt{-1} \partial 
\bar{\partial} u \right).
\end{equation}
As $\widetilde{\omega}_t - \omega_t -  \sqrt{-1} \partial 
\bar{\partial} u$ vanishes at $t=0$, it follows that
$\widetilde{\omega}_t = \omega_t +  \sqrt{-1} \partial \bar{\partial} 
u$ for all $t$. Thus $u$ satisfies (\ref{eq:skrf2}).

Conversely, suppose that we have a smooth solution to 
(\ref{eq:skrf2}) on a time interval $[0,T]$ so that each 
$\omega_t + \sqrt{-1} \partial \bar{\partial} u$ is a K\"ahler 
metric which is biLipschitz equivalent to $\omega_0$, 
and for each $k$, the $k$-th covariant derivatives of $u$ 
(with respect to the initial metric $\omega_0$) are uniformly 
bounded. Putting $\widetilde{\omega}_t = \omega_t + \sqrt{-1} 
\partial \bar{\partial} u$ gives a solution to (\ref{eq:krf}) on
$[0,T]$. Because $\widetilde{\omega}_t$ is biLipschitz equivalent 
to $\omega_0$, each time slice is complete.  From the 
derivative bounds on $u$, the curvature of $\widetilde{\omega}_t$ 
is uniformly bounded on $[0, T]$.
\end{proof}

\begin{remark}
In view of the uniqueness of $\widetilde{\omega}$ on 
$[0, T]$, the uniqueness of $u$ on $[0,T]$ is an 
immediate consequence, since $u$ must satisfy 
the equation (\ref{ueqn2}) and hence be given 
by (\ref{ueqn1}).
\end{remark}

\section{Existence result }

In this section we characterize the first singularity
time for a K\"ahler-Ricci flow solution on a general complex manifold.
Recall the definition of $\omega_t$ from (\ref{omegat}).

\begin{theorem}
\label{th:existence1}
Suppose that $\omega_0$ is a complete K\"ahler metric on a
complex manifold $X$, with bounded curvature.

Let $T_1$ be the supremum (possibly infinite) of the numbers 
$T^\prime$ so that there is a smooth solution for $u$ in 
(\ref{eq:skrf2}) on the time interval $[0,T^\prime]$ such that 
\begin{enumerate}
\item For each $t \in [0, T^\prime]$, $\omega_t + \sqrt{-1} 
\partial \bar{\partial} u$ is a K\"ahler metric which is 
biLipschitz equivalent to $\omega_0$ and
\item For each $k$, the $k$-th covariant derivatives of $u$ 
(with respect to the initial metric $\omega_0$) are uniformly 
bounded on $[0, T^\prime]$.
\end{enumerate}

Let $T_2$ be the supremum (possibly infinite) of the numbers 
$T$ for which there is a function $F_{T} \in C^\infty(X)$ such 
that 
\begin{enumerate} 
\setcounter{enumi}{2}
\item $\omega_{T} + \sqrt{-1} \partial \bar{\partial} F_{T}$ is a
K\"ahler metric which is biLipschitz equivalent to $\omega_0$ 
and
\item For each $k$, the $k$-th covariant derivatives of $F_{T}$ 
(with respect to the initial metric $\omega_0$) are uniformly 
bounded.
\end{enumerate}

Then $T_1 = T_2$.
\end{theorem}
\begin{proof}
If there is a solution for $u$ in (\ref{eq:skrf2}) on a time 
interval $[0,T^\prime]$ satisfying (1) and (2) then we 
can take $F_{T^\prime} = u(T^\prime)$ to show that 
$T_2 \geq T_1$. Thus it suffices to show that $T_1 
\geq T_2$. That is, we need to show that if we can find a 
function $F_T$ satisfying (3) and (4) then we can 
solve $u$ in (\ref{eq:skrf2}) on the time interval $[0, T)$ 
so that for each $T^\prime \in [0, T)$, the restriction of 
the solution to $[0, T^\prime]$ satisfies (1) and (2).

We know that there is a solution $u$ for short time satisfying
(1) and (2). Suppose initially that 
we have a solution on some time interval $[0, \widehat{T})$, with 
$\widehat{T}  < T$, so that (1) and (2) are satisfied on subintervals 
$[0, T^\prime] \subset [0, \widehat{T})$. Our goal is to
derive uniform estimates for the solution $u$ and $\widetilde\omega
_t$, i.e. to show that there are positive numbers $C > 1$ and 
$\{A_k\}_{k=0}^\infty$
so that for all $t \in [0, \widehat{T})$ and $k \ge 0$, 
$\sup_{x \in X} |\nabla^k u|(t,x) 
\leq A_k$ and $C^{-1}\omega_0\leq\widetilde\omega_t
\leq C\omega_0$. 

{\it Note: in what follows, $C$ always stands for a positive constant,
which might be different from place to place.} 

We now give some equations derived from (\ref{eq:skrf2}), as in
\cite{t-znote}. All of the inner products and Laplacians are 
computed using the metric $\widetilde\omega_t = \omega_t + 
\sqrt{-1} \partial \bar{\partial} u$. We also use the fact that
\begin{equation}
\triangle u = \Tr \left( \widetilde\omega_t^{-1} \sqrt{-1} \partial
\bar{\partial} u \right), 
\end{equation}
where $\Delta$ stands for the Laplacian operator with 
respect to the flow metric $\widetilde\omega_t$. 

First, the $t$-derivative of (\ref{eq:skrf2}) gives
\begin{equation}
\label{eq:t-derivative}
\frac{\partial}{\partial t} \left( \frac{\partial u}{\partial t} \right) 
=\Delta\(\frac{\partial u}{\partial t}\)-e^{-t} \Tr \left( \widetilde
\omega_t^{-1} (\omega_0+{\rm Ric}(\omega_0)) \right) 
-\frac{\partial u}{\partial t}.
\end{equation}
This implies that
\begin{equation}
\label{e^tu'} 
\frac{\partial}{\partial t}\(e^t\frac{\partial u}{\partial t}\)=\Delta
\(e^t\frac{\partial u}{\partial t}\)- \Tr \left( \widetilde\omega_t^{-1} 
(\omega_0+{\rm Ric}(\omega_0)) \right).
\end{equation}
Also, since
\begin{align}
n & = \Tr \left(  \widetilde\omega_t^{-1}  \widetilde\omega_t \right) =
\Tr \left(  \widetilde\omega_t^{-1} (\omega_t + \sqrt{-1} \partial
\bar{\partial} u) \right)  \\
&= \Tr \left(  \widetilde\omega_t^{-1} ( - \Ric(\omega_0) + e^{-t} 
(\omega_0 + \Ric(\omega_0) ) \right)+ \triangle u, \notag
\end{align}
we get
\begin{equation}
\label{eq:volume} 
\frac{\partial}{\partial t}\(\frac{\partial u}{\partial t}+u\)=\Delta
\(\frac{\partial u}{\partial t}+u\)-n- \Tr \left( \widetilde\omega_t^{-1} 
{\rm Ric}(\omega_0)\right).
\end{equation}
A linear combination of (\ref{e^tu'}) and (\ref{eq:volume}) 
gives that for any $T > 0$,  
\begin{equation} 
\label{eq:finite-time}
\frac{\partial}{\partial t}\((1-e^{t-T})\frac{\partial u}{\partial t}+u
\)=\Delta\((1-e^{t-T})\frac{\partial u}{\partial t}+u\)-n+ \Tr \left( 
\widetilde\omega_t^{-1}  \omega_T\right).
\end{equation}
(Equation (\ref{eq:volume}) can be viewed as the limiting case 
of equation (\ref{eq:finite-time}) when
$T \rightarrow \infty$.)

Next, the $t$-derivative of (\ref{eq:t-derivative}) gives 
\begin{align}
\frac{\partial}{\partial t}\(\frac{\partial^2 u}{\partial t^2}\)= 
& \Delta\(\frac{\partial^2 u}{\partial t^2}\)- \Tr \left( \widetilde
\omega_t^{-1} \frac{\partial \widetilde\omega_t}{\partial t} 
\widetilde\omega_t^{-1} \sqrt{-1} \partial \bar{\partial} 
\frac{\partial u}{\partial t}
\right) + \\
& +e^{-t} \Tr \left( \widetilde\omega_t^{-1} (\omega_0+{\rm Ric}
(\omega_0)) \right) + e^{-t} \Tr \left( \widetilde\omega_t^{-1}\frac
{\partial \widetilde\omega_t}{\partial t}\widetilde\omega_t^{-1} 
(\omega_0+{\rm Ric}(\omega_0)) \right)
+ \notag \\
& -\frac{\partial^2 u}{\partial t^2}. \notag
\end{align}
As 
\begin{equation}
- \Tr \left( \widetilde\omega_t^{-1} \frac{\partial \widetilde\omega_t}
{\partial t} \widetilde\omega_t^{-1} \left[ \sqrt{-1} \partial \bar
{\partial} \frac{\partial u}{\partial t} -e^{-t} (\omega_0+
{\rm Ric}(\omega_0)) \right] \right) = - \Tr \left( \widetilde
\omega_t^{-1} \frac{\partial \widetilde\omega_t}{\partial t} 
\right)^2 =- \left|\frac{\partial\widetilde\omega_t}{\partial t}
\right|^2,
\end{equation}
we obtain
\begin{equation}
\frac{\partial}{\partial t}\(\frac{\partial^2 u}{\partial t^2}\)= 
\Delta\(\frac{\partial^2 u}{\partial t^2}\)+ e^{-t} \Tr \left( 
\widetilde\omega_t^{-1} (\omega_0+{\rm Ric}(\omega_0)) 
\right)- \frac{\partial^2 u}{\partial t^2} - 
\left|\frac{\partial\widetilde\omega_t}{\partial t}\right|^2.
\end{equation}
Its summation with (\ref{eq:t-derivative}) is  
\begin{equation}
\label{eq:vol-decreasing}
\frac{\partial}{\partial t}\(\frac{\partial^2 u}{\partial t^2}+
\frac{\partial u}{\partial t}\)=\Delta\(\frac{\partial^2 u}
{\partial t^2}+\frac{\partial u}{\partial t}\)-\(\frac{\partial^2 
u}{\partial t^2}+\frac{\partial u}{\partial t}\)-\left|\frac
{\partial\widetilde\omega_t}{\partial t}\right|^2.
\end{equation}

By assumption the curvature is uniformly bounded on 
compact intervals of $[0, \widehat{T})$. Hence we can 
apply the maximum principle freely on such intervals. 
Applying it to (\ref{eq:vol-decreasing}) gives
\begin{equation} \label{voldec2}
\frac{\partial^2 u}{\partial t^2}+\frac{\partial u}{\partial t}
\leq Ce^{-t},
\end{equation}
where $C = \sup_{x \in X} \left( \frac{\partial^2 u}
{\partial t^2}+\frac{\partial u}{\partial t} \right)(0,x)$.
From (\ref{ueqn2}), 
\begin{equation} \label{voleqn}
\widetilde\omega^n_t=e^{\frac{\partial u}{\partial t}+u}
\omega_0^n.
\end{equation}
Hence (\ref{voldec2}) indicates the ``essential decreasing'' 
of the volume form along the flow, i.e.    
\begin{equation}
\frac{\partial}{\partial t}\(\frac{\partial u}{\partial t}+u\)
\leq Ce^{-t}.
\end{equation}
Equivalently, $\frac{\partial}{\partial t}\(e^t\frac{\partial u}
{\partial t}\)\leq C$, so
\begin{equation} 
\label{eq:u'-upper}
\frac{\partial u}{\partial t}\leq Cte^{-t},
\end{equation}
which implies that
\begin{equation} \label{uupper}
u\leq C.
\end{equation} 

To get a lower bound on $u$, we use (\ref{eq:finite-time}). 
We have a smooth bounded function $F_T$ so that
$\omega_T+\sqrt{-1}\partial\bar\partial F_T$ is a K\"ahler 
metric. Then (\ref{eq:finite-time}) can be reformulated as 
\begin{align} 
\frac{\partial}{\partial t}\((1-e^{t-T})\frac{\partial u}{\partial t}+u
-F_T\)= 
& \Delta\((1-e^{t-T})\frac{\partial u}{\partial t}+u-F_T\)-n+ \\
& +\Tr \left( \widetilde\omega_t^{-1} ( \omega_T+\sqrt{-1} \p\bar\p 
F_T ) \right). 
\notag
\end{align}
The maximum principle gives  
\begin{equation} \label{inter}
(1-e^{t-T})\frac{\partial u}{\partial t}+u - F_T + nt \geq 
- \sup F_T.
\end{equation}
Equations (\ref{eq:u'-upper}) and (\ref{inter}) imply a uniform
lower bound 
for $u$ on $[0, \widehat{T})$, of the form
\begin{equation}
u \ge  - C t e^{-t} (1-e^{t-T}) + F_T - \sup F_T - nt.
\end{equation}
Also, equations
(\ref{uupper}) and (\ref{inter}) imply a uniform lower bound for 
$\frac{\partial u}{\partial t}$ on $[0, \widehat{T})$, of the form
\begin{equation}
\frac{\partial u}{\partial t} \ge 
\frac{-C + F_T - \sup F_T - nt}{1-e^{t-T}}.
\end{equation}  

So far we have obtained $0$-th order estimates for
$u$ and $\frac{\partial u}{\partial t}$. 
In order to get uniform higher order estimates (for $u$) 
on $[0, \widehat{T})$, we modify the background form $\omega_t$ 
to make it a K\"ahler form.

First, for any $t \in [0, \widehat{T})$, one has
\begin{equation} \label{put}
\omega_t=\left(\frac{e^{-t}-e^{-T}}{1-e^{-T}}\right)
\omega_0+\left(\frac{1-e^{-t}}{1-e^{-T}}\right)\omega
_T.
\end{equation}
Putting
\begin{equation} \label{put2}
\widehat\omega_t=\left(\frac{e^{-t}-e^{-T}}{1-e^{-T}}\right)\omega
_0+\left(\frac{1-e^{-t}}{1-e^{-T}}\right)(\omega_T+\sqrt{-1}
\partial\bar\partial F_T)
\end{equation}
gives a K\"ahler form which, in fact, is uniformly K\"ahler for
$t \in [0, \widehat{T})$. If we put
\begin{equation} \label{put3}
v=u-\left(\frac{1-e^{-t}}{1-e^{-T}}\right)F_T
\end{equation}
then from (\ref{put}), (\ref{put2}) and (\ref{put3}),
\begin{equation} \label{put4}
\widetilde{\omega}_t = \omega_t+\sqrt{-1}\partial\bar\partial u=\widehat
\omega_t+\sqrt{-1}\partial\bar\partial v.
\end{equation}
The evolution equation 
for $v$ is 
\begin{equation}
\label{eq:skrf3}
\frac{\partial v}{\partial t}={\rm log}\frac{(\widehat\omega_t+\sqrt{-1}
\partial\bar\partial v)^n}{\omega^n_0}-v-\frac{F_T}{1-e^{-T}},
\: \: \: \: v(0, \cdot)=0.  
\end{equation}
The initial value for $\frac{\partial v}{\partial t}$ is no longer 
zero, but this will not cause any problems. The extra term 
$\frac{F_T}{1-e^{T}}$ in the right-hand side of (\ref{eq:skrf3}) 
is also well controlled. 
The point is that
equation (\ref{eq:skrf3}) is phrased 
in terms of a background metric $\widehat\omega_t$ which is uniformly 
K\"ahler on $[0, \widehat{T})$. The $0$-th order bounds on $u$ 
and $\frac{\partial u}{\partial t}$ imply $0$-th order bounds on 
$v$ and $\frac{\partial v}{\partial t}$. \\

In the following, we sketch the argument to obtain the 
higher order estimates, which is fairly standard.

We begin with an estimate on $\triangle_{\widehat\omega} v$.
We use computations
in \cite{yau} to derive an inequality for solutions of
(\ref{eq:skrf3}). This inequality is closely related to
to \cite[(1.5)]{cao}, \cite[(2.3)]{t-znote}) and \cite[(15)]{tsu}.
\begin{lemma} \label{2.24}
Given a solution of (\ref{eq:skrf3}),
there is an inequality of the form
\begin{align}
e^{Cv}\(\Delta-\frac{\partial}{\partial{t}}\)\bigl( e^{-Cv}
(n+\Delta_{\widehat\omega_t} v)\bigr) \geq 
& \Delta_{\widehat\omega_{t}}\({\rm log}\frac{
\Omega}{\widehat\omega^n_t}\)-n^{2} \inf_{i\neq j}R_{i\bar{i}
j\bar{j},t}-n + \\
& \(C\frac{\partial{v}}{\partial{t}}-C\)
(n+\Delta_{\widehat\omega} v)+ \notag \\
& (C+ \inf_{i\neq j}R_{i\bar{i}j\bar{j},t})e^{-\frac{\frac
{\partial{v}}{\partial {t}}+v+{\rm log}\frac{\Omega}
{\widehat\omega^n_t}}{n-1}}{(n+\Delta_{\widehat\omega} v)}^
{\frac{n}{n-1}}, \notag
\end{align}
where 
\begin{equation} \label{added}
n+\Delta_{\widehat\omega_t} v=\Tr \(\widehat\omega
^{-1}_t(\widehat\omega_t+\sqrt{-1}\partial\bar\partial v)\)
=\Tr \(\widehat\omega^{-1}_t\widetilde\omega_t\)>0,
\end{equation}
\begin{equation}
\Omega=e^{\frac{F_T}{1-e^{-T}}}\omega^n_0
\end{equation}
and $C$ is a constant that depends on $t$, $\omega_0$ 
and $F_T$.
\end{lemma}
\begin{proof}
As in \cite[Section 2]{yau}, suppose that $\phi$ is a smooth 
solution of
\begin{equation}
(\omega+\sqrt{-1}\partial\bar{\partial}\phi)^n=e^{f}
{\omega}^n,
\end{equation}
where $\omega$ is a K\"ahler metric on a K\"ahler manifold 
$X$ and $f$ is a smooth function. We have the following 
inequality at any point $p \in X$:
\begin{align} \label{2.29}
e^{C\phi}\Delta\bigl( e^{-C\phi}(n+\Delta
_\omega\phi)\bigr) \geq & \Delta_{\omega}f-n^{2}
\inf_{i\neq j}R_{i\bar{i}j\bar{j}}-Cn(n+\Delta_\omega\phi)+ \\
& (C+\inf_{i\neq j}R_{i\bar{i}j\bar{j}})e^{-\frac{F}{n-1}}{(n+
\Delta_\omega\phi)}^{\frac{n}{n-1}}. \notag
\end{align}
Here $R_{i\bar{i}j\bar{j}}$ comes from the curvature tensor 
for the metric $\omega$ (written in terms of any unitary 
frame $\{e_i\}_{i=1}^n$), $C$ is any (fixed) positive 
constant such that $C+\inf_{i\neq j}R_{i\bar{i}j\bar{j}}>0$ 
at $p$, and $\Delta$ is the Laplacian with respect to 
$\omega+\sqrt{-1}\partial\bar\partial \phi$.
We emphasize that this inequality is pointwise and 
the ``$\inf$'' is taken at the point $p$. 

Now we consider the flow. Equation (\ref{eq:skrf3}) 
can be reformulated as
\begin{equation} \label{2.30}
(\widehat\omega_{t}+\sqrt{-1}\partial\bar{\partial}v)^n=e^{\frac
{\partial{v}}{\partial{t}}+v+{\rm log}\frac{\Omega}{\widehat\omega
^n_t}}\widehat\omega^n_t.
\end{equation}
From (\ref{2.29}),
\begin{align} \label{2.31}
e^{Cv}\Delta\bigl( e^{-Cv}(n+\Delta_{\widehat\omega_t} v)\bigr) 
\geq & \Delta_{\widehat\omega_{t}} \left( \frac{\partial{v}}
{\partial{t}}+v+{\rm log}\frac{\Omega}{\widehat\omega^n_t} \right) -n^{2}
\inf_{i\neq j}R_{i\bar{i}j\bar{j}, t}  - \\
& Cn(n+\Delta_{\widehat\omega_t} v)+ \notag \\
& (C+\inf_{i\neq j}R_{i\bar{
i}j\bar{j}})e^{-\frac{\frac{\partial{v}}{\partial{t}}+v+{\rm log}
\frac{\Omega}{\widehat\omega^n_t}}{n-1}}{(n+\Delta_{\widehat
\omega_t} v)}^{\frac{n}{n-1}}, \notag
\end{align}
where $R_{i\bar{i}j\bar{j}, t}$ is computed using the metric
$\widehat\omega_{t}$, $C$ is a positive constant such that 
$C+\inf_{i\neq j}R_{i\bar{i}j\bar{j}, t}>0$
and $\Delta$ is  the Laplaciain with respect to the flow metric $\widehat\omega_t+\sqrt{-1}\partial\bar\partial v$.
Next, one computes that
\begin{align} \label{2.32}
e^{Cv} \left( -\frac{\partial}{\partial{t}} \right) \left( e^{-Cv}(n+
\Delta_{\widehat\omega_t} v) \right) 
=& C\frac{\partial{v}}{\partial{t}}(n+\Delta_{\widehat\omega_t} v)-
\frac{\partial}{\partial{t}}(n+\Delta_{\widehat\omega_t} v) \\
=& C\frac{\partial{v}}{\partial{t}}(n+\Delta_{\widehat\omega_t} v)-
\frac{\partial}{\partial{t}} \Tr \(
\widehat\omega^{-1}_{t} \sqrt{-1}\partial\bar{\partial} v\) \notag \\
=& C\frac{\partial{v}}{\partial{t}}(n+\Delta_{\widehat\omega_t} v)+
 \Tr \( \widehat\omega^{-1}_{t} \frac{\partial \widehat\omega_{t}}{\partial t}\widehat\omega^{-1}_{t} \sqrt{-1}\partial\bar{\partial} v\) -\notag \\
&\Delta_{\widehat\omega
_{t}}(\frac{\partial {v}}{\partial {t}}).\notag
\end{align}

Adding (\ref{2.31}) and (\ref{2.32}) gives 
\begin{align}
e^{Cv} \left( \Delta-\frac{\partial}{\partial{t}} 
\right) \bigl( e^{-Cv}(n+\Delta_{\widehat\omega_t} v)\bigr) 
\geq & \Delta_{\widehat\omega_{t}}\({\rm log}\frac{
\Omega}{\widehat\omega^n_t}\)-n^{2} \inf_{i\neq j}R_{i\bar{i}
j\bar{j}, t}-n + \\
& \(-Cn+C\frac{\partial{v}}{\partial{t}}+1\)
(n+\Delta_{\widehat\omega_t} v)+ \notag \\
&\Tr \( \widehat\omega^{-1}_{t} \frac{\partial \widehat\omega_{t}}{\partial t}
\widehat\omega^{-1}_{t} \sqrt{-1}\partial\bar{\partial} v\) + \notag \\
&(C+\inf_
{i\neq j}R_{i\bar{i}j\bar{j}, t})e^{-\frac{\frac{\partial{v}}{
\partial {t}}+v+{\rm log}\frac{\Omega}{\widehat\omega^n_t}}
{n-1}}{(n+\Delta_{\widehat\omega_t} v)}^{\frac{n}{n-1}}. \notag
\end{align}
Since $\frac{\partial{\widehat\omega_{t}}}{\partial{t}}$ is relatively bounded with respect to $\widehat\omega_t$, 
equation (\ref{added}) implies that
\begin{equation}
\Tr \( \widehat\omega^{-1}_{t} \frac{\partial \widehat\omega_{t}}{\partial t}
\widehat\omega^{-1}_{t} \sqrt{-1}\partial\bar{\partial} v\) \geq -C(n+\Delta_{\widehat\omega} v)-C
\end{equation}
for some $C>0$.
Hence after a redefinition of $C$, we have
\begin{align}
e^{Cv} \left( \Delta-\frac{\partial}{\partial{t}} \right) \bigl( e^{-Cv}
(n+\Delta_{\widehat\omega} v)\bigr) 
\geq &\Delta_{\widehat\omega_{t}}\({\rm log}\frac{
\Omega}{\widehat\omega^n_t}\)-n^{2} \inf_{i\neq j}R_{i\bar{i}
j\bar{j}, t}- C +\\
& \(C\frac{\partial{v}}{\partial{t}}-C\)
(n+\Delta_{\widehat\omega} v)+\notag \\
&(C+ \inf_{i\neq j}R_{i\bar{i}j\bar{j}, t})e^{-\frac{\frac
{\partial{v}}{\partial {t}}+v+{\rm log}\frac{\Omega}
{\widehat\omega^n_t}}{n-1}}{(n+\Delta_{\widehat\omega} v)}
^{\frac{n}{n-1}}. \notag
\end{align}
This proves the lemma.
\end{proof}

Using Lemma \ref{2.24} and the $0$-th order bounds on $v$ 
and $\frac{\partial v}{\partial t}$, along with the uniform control 
on $\widehat\omega_t$ as a metric, we conclude that there is an 
estimate of the form
\begin{align} \label{discussion}
\left( \Delta-\frac{\partial}{\partial {t}} \right) \bigl( e^{-Cv}(n+
\Delta_{\widehat\omega_t} v)\bigr) 
&\geq -C+\(C\frac{\partial {v}}{\partial {t}}-C\)(n+
\Delta_{\widehat\omega_t} v)+C(n+\Delta_{\widehat\omega_t} 
v)^{\frac{n}{n-1}} \\
&\geq -C-C(n+\Delta_{\widehat\omega_t} v)+C(n+
\Delta_{\widehat\omega_t} v)^{\frac{n}{n-1}}. \notag
\end{align}
The maximum principle now
gives an {\it a priori} upper bound for $e^{-Cv}
(n+\Delta_{\widehat\omega_t} v)$, and hence also for $(n+\Delta_
{\widehat\omega_t} v)$. 

This Laplacian upper bound gives a trace upper bound on $\widetilde\omega_t$, relative
to $\widehat\omega_t$, from
(\ref{put4}). There is also a determinant lower bound on  $\widetilde\omega_t$ relative
to $\widehat\omega_t$, coming from (\ref{2.30}), where we use the lower bounds on
$v$ and $\frac{\p v}{\p t}$.  This gives a uniform bound on $\widetilde\omega_t$, relative
to $\widehat\omega_t$.

Next, we look at  third order estimates. Following \cite{cao} and
\cite{yau}, we consider the expression $S=
\widetilde g^{i\bar j}\widetilde g^{k\bar l}\widetilde g^{\lambda\bar
\eta}v_{i\bar l\lambda}v_{\bar jk\bar\eta}$, where 
$\widetilde g_{i\bar j}$ is the metric tensor corresponding to $\widetilde
\omega_t$. As in 
\cite[Section 5.3]{chau}, there are estimates of the form
\begin{align}
\(\Delta-\frac{\partial}{\partial t}\)S 
&\geq -C\cdot S-C, \\
\(\Delta-\frac{\partial}{\partial t}\)\Delta_{\widehat
\omega_{t}} v 
&\geq  C\cdot S-C. \notag
\end{align}

Choosing $A>0$ large enough, as in
\cite[(1.25)]{cao} and \cite[(5.17)]{chau}, there is an estimate of the form
\begin{equation}
\(\Delta-\frac{\partial}{\partial t}\)(S+A\Delta_{\widehat
\omega_{t}} v)\geq C\cdot S-C.
\end{equation}
Applying the maximum principle and using the uniform 
control on $\Delta_{\widehat\omega_t} v$, we obtain an {\it a priori} 
upper bound on $S$.
This provides a spatial $C^{2, \alpha}$-bound for $v$ and a
$C^\alpha$-bound for the metric coefficients of $\widetilde
\omega_t$. 

One can then obtain further derivative bounds 
(cf. \cite[Section 5.4]{chau}) and apply parabolic 
Schauder estimates (cf. \cite[Section 5.5]{chau}). 
In this way, one obtains {\it a priori} estimates on 
all of the derivatives of $v$. These imply the 
desired derivative estimates on $u$.

To summarize, we assumed that we have a solution of 
(\ref{eq:skrf2}) on a time interval $[0, \widehat{T})$, with 
$\widehat{T} < T$, so that conditions (1) and (2) are 
satisfied on compact subintervals of $[0, \widehat{T})$. 
Then we have shown that there are numbers 
$C> 1$ and $\{A_k\}_{k=0}^\infty$ so that for all $t \in [0, \widehat
{T})$ and ${x \in X}$, we have $|\nabla^k u|(t,x) \le A_k$ 
and $C^{-1}\omega_0\leq\widetilde\omega_t\leq C\omega_0$. 

Let $\{t_i\}_{i=1}^\infty$ be a sequence in 
$[0, \widehat{T})$ with $\lim_{i \rightarrow \infty} t_i = 
\widehat{T}$. We can extract a subsequence of $\{u(t_i, 
\cdot)\}_{i=1}^\infty$ that converges in the pointed 
$C^\infty$-topology to some $u_{\widehat{T}}(\cdot) \in C^\infty
(X)$. Now $u_{\widehat{T}}$ is uniformly bounded on $X$, along 
with its covariant derivatives (with respect to $\omega_0$). 
From (\ref{voleqn}), it follows that $\omega_{\widehat{T}} 
+ \sqrt{-1}\partial \bar{\partial} u_\infty$ is K\"ahler and 
biLipschitz to $\omega_0$, with bounded curvature. 
Hence we can solve the equation in (\ref{eq:skrf2}) to get a solution 
$U$ on a time interval $[\widehat{T}, \widehat{T} + 
\epsilon)$ with initial condition $U(\widehat{T}) = u_{\widehat{T}}$. 
One then shows that the solutions $u(\cdot)$, on $[0, 
\widehat{T})$, and $U(\cdot)$, on $[\widehat{T}, \widehat
{T} + \epsilon)$, join to form a smooth solution of 
(\ref{eq:skrf2}) on $[0, \widehat{T} + \epsilon)$, which 
satisfies conditions (1) and (2) on compact subintervals. 
It follows that there is a solution of (\ref{eq:skrf2}) on 
$[0, T)$ so that for each $T^\prime \in [0, T)$, the 
restriction of the solution to $[0, T^\prime]$ satisfies 
conditions (1) and (2). This finishes the proof.
\end{proof}

\section{Long-time convergence}

In this section we prove a long-time convergence result for
the normalized K\"ahler-Ricci flow equation,
under the assumption that the initial metric satisfies
an inequality of the form
$-{\rm Ric}(\omega_0)+\sqrt{-1}\partial\bar\partial 
f>\epsilon\omega_0$ for some bounded function $f$ and some
positive constant $\epsilon$.
We show that the solution smoothly 
approaches a complete K\"ahler-Einstein metric,
having Einstein constant $-1$.
This result can be seen as a
generalization of results in \cite{cao}, \cite{t-znote} concerning
K\"ahler-Ricci flow solutions on compact manifolds. It is also
a generalization of \cite[Theorem 1.1]{chau}, 
which proves the same conclusion under the assumption
that $-{\rm Ric}(\omega_0)+\sqrt{-1}\partial\bar\partial 
f=\omega_0$.

\begin{theorem}
\label{th:non-deg-convergence1}
1. Suppose that $\omega_0$ is a complete K\"ahler metric on a
complex manifold $X$, with bounded curvature, such that
$-{\rm Ric}(\omega_0)+\sqrt{-1}\partial\bar\partial 
f \ge 0$ for some smooth function $f$ with bounded 
$k$-th covariant derivatives (with respect to $\omega_0$) for
each $k \ge 0$. Then the flow (\ref{eq:krf}) (or equivalently 
(\ref{eq:skrf2})) exists forever.  \\
2. Suppose in addition that
$-{\rm Ric}(\omega_0)+\sqrt{-1}\partial\bar\partial 
f>\epsilon\omega_0$ for some $\epsilon > 0$.
Then the flow (\ref{eq:krf}) (or equivalently 
(\ref{eq:skrf2}))
converges smoothly to a complete K\"ahler-Einstein metric 
with Einstein constant $-1$. 
\end{theorem}

\begin{proof}
Suppose first that $-{\rm Ric}(\omega_0)+\sqrt{-1}\partial\bar\partial 
f \ge 0$.
Then
\begin{equation}
\omega_t + \left( 1 - e^{-t} \right) \sqrt{-1} \partial
\overline{\partial} f = 
e^{-t} \omega_0 + \left( 1 - e^{-t} \right) \left(
-\Ric(\omega_0) + \sqrt{-1} \partial
\overline{\partial} f \right) \ge e^{-t} \omega_0.
\end{equation}
From Theorem \ref{th:existence1}, the flow (\ref{eq:krf}) exists forever. 

Now suppose that 
$-{\rm Ric}(\omega_0)+\sqrt{-1}\partial\bar\partial 
f>\epsilon\omega_0$.
To prove the long-time convergence,
we need estimates that are uniform in time. The upper 
bounds on $u$ and $\frac{\partial u}{\partial t}$ 
from (\ref{uupper}) and (\ref{eq:u'-upper})
are uniform 
for all time. For the lower bound, we use the following 
variation on (\ref{eq:volume}): 
\begin{equation} \label{variation}
\frac{\partial}{\partial t}\(\frac{\partial u}{\partial t}+u-f\)=\Delta
\(\frac{\partial u}{\partial t}+u-f\)-n+\Tr\(\widetilde\omega^{-1}_t
(-{\rm Ric}(\omega_0)+\sqrt{-1}\p\bar\p f)\). 
\end{equation}
Now
\begin{align}
\Tr\(\widetilde\omega^{-1}_t(-{\rm Ric}(\omega_0)+\sqrt{-1}\p\bar\p 
f)\)
&\geq n\cdot \left(\frac{(-{\rm Ric}(\omega_0)+\sqrt{-1}\p
\bar\p f)^n}{\widetilde\omega_t^n}\right)^{\frac{1}{n}} \\
&\geq n\cdot \left(\frac{(\epsilon\omega_0)^n}{\widetilde
\omega_t^n}\right)^{\frac{1}{n}} \notag \\
&= n \epsilon \: e^{-\frac{1}{n}(\frac{\partial u}{\partial t}+u)}, \notag
\end{align}
so (\ref{variation}) gives
\begin{equation} \label{variation2}
\frac{\partial}{\partial t}\(\frac{\partial u}{\partial t}+u-f\) \ge \Delta
\(\frac{\partial u}{\partial t}+u-f\)-n+n \epsilon \: 
e^{- \frac{f}{n}} \:
e^{-\frac{1}{n}(\frac{\partial u}{\partial t}+u-f)}.
\end{equation}
Putting $Y(t) = \inf_{x \in X} \left(\frac{\partial u}{\partial t}+u-f
\right)(x,t)$, we can apply the maximum principle to (\ref{variation2}) to
conclude that $Y(t)$ is bounded below by the solution $c(t)$ to the
ordinary differential equation 
\begin{equation}
\frac{dc}{dt} = -n + n \epsilon \: e^{-\frac{\sup f}{n}} \: e^{- \frac{c}{n}}
\end{equation}
with initial condition $c(0) = Y(0)$.
It follows that there is a lower bound
\begin{equation} \label{lowerb}
\frac{\partial u}{\partial t}+u\geq -C
\end{equation}
which is uniform in $t$. When combined with the upper bounds on
$u$ and $\frac
{\partial u}{\partial t}$, equation (\ref{lowerb}) 
provides uniform lower bounds for $u$ and $\frac
{\partial u}{\partial t}$.

As in the proof of Theorem \ref{th:existence1}, we now transform 
the flow equation in order to prove the 
higher order estimates.
Putting
\begin{equation} \label{backgroundomegat}
\widehat\omega_t=
e^{-t} \omega_0 +
(1-e^{-t}) (-{\rm Ric}(\omega_0)+\sqrt{-1}\partial\bar
\partial f),
\end{equation}
we have 
\begin{equation}
\widehat\omega_t \ge \left( e^{-t} + \epsilon (1 - e^{-t}) \right) 
\omega_0,
\end{equation}
so the family $\{\widehat\omega_t\}_{t \ge 0}$ is uniformly K\"ahler.
Next, putting
\begin{equation}
w=u-(1-e^{-t})f,
\end{equation}
we can write
\begin{equation}
\widetilde\omega_t=\omega_t+\sqrt{-1}\partial\bar\partial 
u=\widehat\omega_t+\sqrt{-1}\partial\bar\partial 
w.
\end{equation}
Then the flow equation (\ref{eq:skrf2}) becomes
\begin{equation}
\label{eq:skrf4}
\frac{\partial w}{\partial t}={\rm log}
\frac{(\widehat\omega_t+\sqrt{-1}\partial\bar\partial w)^n}
{\omega^n_0}-w-f, ~~~~ w(0, \cdot)=0.  
\end{equation}
We can use this equation to find higher order estimates on $u$, as in the
proof of Theorem \ref{th:existence1}.
Note that the background metric $\widehat\omega_t$ is uniformly
K\"ahler and from (\ref{backgroundomegat}), it is uniformly bounded above. 
Hence
the higher order estimates will also be uniform in
time. So we have achieved uniform estimates on $u$ for all 
time.

We now justify the convergence. Using the uniform bounds that
we have obtained so far,
equation (\ref{eq:t-derivative}) implies an inequality of the form
\begin{equation}
\frac{\partial}{\partial t} \left( e^t \frac{\partial u}{\partial t} \right) 
\ge \Delta\(e^t \frac{\partial u}{\partial t}\)-C.
\end{equation}
From the maximum principle,
$e^t\frac{\partial u}{\partial t}+Ct\geq 0$ and so 
\begin{equation} \label{combine}
\frac{\partial u}{\partial t}\geq -Cte^{-t}.
\end{equation}  
Combining (\ref{eq:u'-upper}) and (\ref{combine}), we conclude that
$\lim_{t \rightarrow \infty} u(x,t) = u_\infty(x)$ for some
function $u_\infty$ on $X$. 
Using the uniform higher order 
bounds on $u(t)$, 
one sees that there is uniform $C^K$-convergence of $u(t)$ toward
$u_\infty$ for any $K > 0$.
Taking the limit of 
(\ref{eq:krf}) 
as $t \rightarrow \infty$ shows that the limiting K\"ahler metric
$\omega_\infty = -\Ric(\omega_0) 
+ \sqrt{-1} \partial \overline{\partial} u_\infty$
satisfies ${\rm Ric}(\omega_\infty)=-\omega_\infty$. 
Also, $\omega_\infty$ is biLipschitz equivalent to
$\omega_0$; see the discussion after (\ref{discussion}) and note that
\begin{equation}
\omega_\infty = (- \Ric(\omega_0) + \sqrt{-1} \partial \overline{\partial} f)
+ \sqrt{-1} \partial \overline{\partial} (u_\infty - f).
\end{equation}
Thus $\omega_\infty$ is complete.
\end{proof}

For an example of Theorem \ref{th:non-deg-convergence1}, see 
Example \ref{ke}.

\section{Standard spatial asymptotics}

In this section we begin to specify the spatial asymptotics
that we want to consider.  The goal is to come up with the
widest class of spatial asymptotics which is preserved by
the Ricci flow, and for which we can prove something nontrivial.
We introduce the notion of ``standard''
spatial asymptotics for a K\"ahler metric
$\omega_X$ on a quasiprojective
manifold $X = \overline{X} - D$.  Assuming standard spatial
asymptotics,
we prove some properties of
the extension of $\omega_X$ by zero to $\overline{X}$. 
As an example of standard spatial asymptotics,
in the case when $[K_{\overline{X}} + D] > 0$, we show how to
recover the K\"ahler-Einstein metric on $X$
\cite{kobayashi,tian-yau,tsu} using the K\"ahler-Ricci flow,
via Theorem \ref{th:non-deg-convergence1} or
\cite[Theorem 1.1]{chau}.

In the next section we will show that the property of having
standard spatial asymptotics is preserved by Ricci flow.  In
Section \ref{super} we will consider a refinement called
``superstandard'' spatial asymptotics.

Suppose that $\overline{X}$ is a compact $n$-dimensional
complex manifold and $D$ is an effective
divisor with simple normal crossings.  Put $X = \overline{X} - D$.
Let $D = \sum_{i=1}^k D_i$ be the decomposition of $D$ into
its irreducible components.  Let $L_i$ be the holomorphic line
bundle on $\overline{X}$ corresponding to $D_i$. 
Put $L_D = \bigotimes_{i=1}^k L_i$. 

Let $h_{L_i}$ be a Hermitian metric on $L_i$.
There is a holomorphic section $\sigma_i$ of $L_i$ whose zero-set is $D_i$,
unique up to multiplication by a nonzero complex number.
The section $\sigma_i$ is nondegenerate along $D_i$
\cite[Theorem II.(6.6)]{demailly}. That is,
the restriction of the bundle map $\nabla^{L_i} \sigma_i : T \overline{X} 
\rightarrow L_i$ to $D_i$ factors through an isomorphism 
$T\overline{X}/TD_i \rightarrow L_i\big|_{D_i}$.

Given a multi-index $I = (i_1, \ldots, i_m)$, put
$D_I = \bigcap_{j=1}^m D_{i_j}$. We write $|I| = m$. Put
$D_I^{int} = D_I - \bigcup_{I' : |I^\prime| > m} D_{I^\prime}$.
Then $D_I^{int}$ is a smooth complex manifold of complex dimension $n-m$,
possibly noncompact.

Let $\Delta$ denote the open unit ball in $\C$. Put
$\Delta^* = \Delta - \{0\}$ and $H = \{ z \in \C : 
\text{Im}(z) > 0\}$.
There is a holomorphic covering map $\pi :  H \rightarrow \Delta^*$
given by $\pi(z) = e^{iz}$.
Suppose that $\overline{x} \in D_I^{int}$. After permutation of indices, we
can assume that 
$\overline{x} \in (D_1 \cap D_2 \cap \ldots \cap D_m) - (D_{m+1} \cup \
D_{m+2} \cup \ldots \cup D_k)$.
We write $0$ for $(0,\ldots,0) \in \Delta^n$.
 Then there is a neighborhood $U$
of  $\overline{x}$ in $\overline{X}$
and a biholomorphic map $F_{\overline{x}} : \Delta^n \rightarrow U$
so that 
\begin{enumerate}
\item  For $i > m$, $U \cap D_i = \emptyset$.
\item $F_{\overline{x}}(0) = \overline{x}$.
\item For $1 \le i \le m$, 
$F_{\overline{x}}(\Delta^{i-1} \times \{0\} \times  \Delta^{n-i}) = U \cap D_i$.
\item For $1 \le i \le m$, 
$\parallel \sigma_i(F_{\overline{x}}(z)) \parallel^2_{L_i} = h_i |z^i|^2$ for some
positive function $h_i \in C^\infty(\Delta^n)$.
\end{enumerate}
In particular, $F_{\overline{x}}  \left( (\Delta^*)^m \times \Delta^{n-m} \right) = U \cap X$.
Passing to the universal cover gives a holomorphic covering map
$\widetilde{F}_{\overline{x}} : H^m \times \Delta^{n-m} \rightarrow U \cap X$.

The map $G_{\overline{x}}$ on $\Delta^{n-m}$, given by
$G_{\overline{x}}( w) = F_{\overline{x}}(0,w)$, is a biholomorphic map from
$\Delta^{n-m}$ to a neighborhood of $\overline{x}$ in 
$D_I^{int}$.

Let $\omega_X$ be a
K\"ahler metric on $X$. Then
$\widetilde{F}_{\overline{x}}^* \omega_X$ is a K\"ahler metric on
$H^m \times \Delta^{n-m}$ which is invariant under translation in the $H^m$-factor
by $2\pi  \Z^m$. Given
$r  \in (\R^+)^m$, define a biholomorphic map $\alpha_r : H^m \rightarrow H^m$ by
$\alpha_r(h_1, \ldots, h_m) = (r_1 h_1, \ldots, r_m h_m)$.
If $Z$ is an auxiliary space then we will also write $\alpha_r$ for
$(\alpha_r, \Id) : H^m \times Z \rightarrow H^m \times Z$.

\begin{definition} \label{standard}
Let $\{\omega_{D_I^{int}} \}$ be
complete K\"ahler metrics on $\{{D_I^{int}} \}$.
Let $\{c_i\}_{i=1}^m$ be positive numbers. Then $\omega_X$ has
standard spatial asymptotics associated to $\{\omega_{D_I^{int}} \}$ and
$\{c_i\}_{i=1}^m$ if for every $\overline{x} \in \overline{X}$
and every local parametrization $F_{\overline{x}}$,
\begin{equation} \label{limitdef}
\lim_{r \rightarrow \infty} \alpha_r^* \widetilde{F}_{\overline{x}}^* \omega_X
\: = \: \sum_{i=1}^m c_i \frac{\sqrt{-1}}{2} \frac{dz^i \wedge
d\overline{z}^i}{(\text{Im}(z^i))^2}  + 
G_{\overline{x}}^* \omega_{D_I^{int}},
\end{equation}
where $\lim_{r \rightarrow \infty}$ means that $r_i \rightarrow \infty$ for each 
$1 \le i \le m$.
The limit in (\ref{limitdef}) is taken in the pointed $C^\infty$-topology around the
basepoint $(\sqrt{-1}, \ldots, \sqrt{-1}) \times 0 \in H^m \times \Delta^{n-m}$.
\end{definition}

\begin{definition}
Let $\{u_{D_I^{int}} \}$ be smooth functions on $\{{D_I^{int}} \}$.
Then a function $u_X \in C^\infty(X)$ has
standard spatial asymptotics associated to $\{u_{D_I^{int}} \}$ 
if for every $\overline{x} \in \overline{X}$
and every local parametrization $F_{\overline{x}}$
\begin{equation} \label{ulimitdef}
\lim_{r \rightarrow \infty} \alpha_r^* \widetilde{F}_{\overline{x}}^* u_X
\: = \: 
G_{\overline{x}}^* u_{D_I^{int}},
\end{equation}
where $\lim_{r \rightarrow \infty}$ means that $r_i \rightarrow \infty$ for each 
$1 \le i \le m$.
The limit in (\ref{limitdef}) is taken in the pointed $C^\infty$-topology around the
basepoint $(\sqrt{-1}, \ldots, \sqrt{-1}) \times 0 \in H^m \times \Delta^{n-m}$.
\end{definition}

\begin{remark}
Note that if $U$ is any neighborhood of $(\sqrt{-1}, \ldots, \sqrt{-1}) \in H^m$
then for some $K_U > 0$,
$\bigcup_{r \ge 1} \alpha_r(U)$ contains $\{z \in \C^m : 
\re (z^i) \in [0, 2\pi],
\text{Im}(z^i) \ge K_U
\text{ for } 1 \le i \le m\}$.
\end{remark}

\begin{theorem} \label{extension}
If $\omega_X$ has standard spatial asymptotics then
\begin{enumerate}
\item $\omega_X$ extends by zero to an element $[\omega_X] \in 
\HH^{(1,1)}(\overline{X}; \R)$.
\item For any $0 \le j \le n$, $\omega_X^j$ extends to an element
of $\HH^{(j,j)}(\overline{X}; \R)$ which equals $[\omega_X]^j$.
\item $\Ric(\omega_X)$ extends by zero to an element 
$[\Ric(\omega_X)] \in \HH^{(1,1)}(\overline{X}; \R)$
which equals $-[K_{\overline{X}} + D]$.
\item $[\omega_X]$ lies in the K\"ahler cone of $\overline{X}$. 
\end{enumerate}
\end{theorem}
\begin{proof}
(1). We use results from \cite[Chapter 3]{demailly}. From Definition 
\ref{standard},
the K\"ahler metric $F_{\overline{x}}^* \omega_X$ is uniformly
biLipshitz equivalent to $2 \sqrt{-1} \sum_{i=1}^m c_i \frac{dz^i \wedge
d\overline{z}^i}{|z^i|^2 \log^2(|z^i|^{-2})} + 
G_{\overline{x}}^* \omega_{D_I^{int}}$ on
$(\Delta^*)^m \times \Delta^{n-m}$.
It follows that $\omega_X$, considered as a current on $\overline{X}$,
has locally finite mass in the sense of \cite[Remark III.(1.15)]{demailly}.
From the Skoda-El Mir extension theorem \cite[Theorem III.(2.3)]{demailly},
$\omega_X$ extends to a closed current on $\overline{X}$ and hence a
class $[\omega_X] \in 
\HH^{(1,1)}(\overline{X}; \R)$. \\
(2). The same argument as in part (1) shows that 
$\omega_X^j$ extends to an element
$[\omega_X^j] \in \HH^{(j,j)}(\overline{X}; \R)$. 
The extension of $\omega_X$
to $\overline{X}$ has no singular support.
Also, $\omega_X^j$ is $L^1$ on $\overline{X}$,
so it is plausible that
$[\omega_X^j] = [\omega_X]^j$.

To see this, let $g_{\overline{X}}$ be a K\"ahler metric
on $\overline{X}$.
Let $\triangle^{(j,j)}$ denote the (nonnegative) Hodge Laplacian on 
$\Omega^{(j,j)}(\overline{X})$ associated to the K\"ahler form
$\omega_{\overline{X}}$.
For small $\epsilon > 0$, the Schwartz kernel 
$e^{-\epsilon \triangle^{(j,j)}}(x,y)$ is well approximated by
$(4 \pi \epsilon)^{-n} \: 
e^{- \frac{d_{\overline{X}}(x,y)^2}{4\epsilon}} P_{x,y}$, where
$P_{x,y}$ denotes parallel transport from $\Lambda^{(j,j)}_y$ to
$\Lambda^{(j,j)}_x$ along a minimal geodesic from $y$ to $x$
(In what
follows, we can
assume that $x$ is not in the cut locus of $y$.) 

For any $\epsilon > 0$,
$e^{-\epsilon \triangle^{(1,1)}} \omega_X$ is a smooth closed form on
$\overline{X}$ whose cohomology class equals
$[\omega_X] \in \HH^{(1,1)}(\overline{X}; \R)$.
Hence 
$\left( e^{-\epsilon \triangle^{(1,1)}} \omega_X \right)^j$ represents
$[\omega_X]^j$, while 
$e^{-\epsilon \triangle^{(j,j)}} \omega_X^j$ represents
$[\omega_X^j]$. 
In particular,
for any smooth closed form $\alpha \in \Omega^{(n-j,n-j)}(\overline{X})$,
we have
\begin{equation}
\int_{\overline{X}} [\omega_X^j] \wedge [\alpha] =
\lim_{\epsilon \rightarrow 0} \int_{\overline{X}}
\left( e^{-\epsilon \triangle^{(j,j)}} \omega_X^j \right) \wedge \alpha =
\lim_{\epsilon \rightarrow 0} \int_{\overline{X}}
\omega_X^j \wedge 
\left( e^{-\epsilon \triangle^{(n-j,n-j)}} \alpha \right) =
\int_{X}
\omega_X^j \wedge \alpha.
\end{equation}

On the other hand,
\begin{equation}
\int_{\overline{X}} [\omega_X]^j \wedge [\alpha] =
\lim_{\epsilon \rightarrow 0} \int_{\overline{X}}
\left( e^{-\epsilon \triangle^{(1,1)}} \omega_X \right)^j \wedge \alpha.
\end{equation}
We claim that there is $L^1$-convergence
\begin{equation} \label{l1}
\lim_{\epsilon \rightarrow 0} 
\left( e^{-\epsilon \triangle^{(1,1)}} \omega_X \right)^j =
\omega_X^j.
\end{equation}
If not then there is a sequence $\{O_k\}_{k=1}^\infty$ 
of nonempty open subsets of $\overline{X}$, 
with $\diam(O_k) \le \frac{1}{k}$,
so that for each $k$ we do not have $L^1$-convergence in
(\ref{l1}) on $O_k$. After passing to a subsequence, we can assume that
there is some $\overline{x} \in \overline{X}$ so that
$\lim_{k\rightarrow \infty} \overline{O_k} = \{ \overline{x} \}$ in the
Hausdorff topology.
Choose a biholomorphic map $F_{\overline{x}} : \Delta^n \rightarrow U$
as before with 
$G_{\overline{x}}(\Delta^{n-m}) = 
F_{\overline{x}}(0, \Delta^{n-m}) \subset D_I^{int}$.
The standard asymptotics from Definition \ref{standard} control
$\omega_X$ on $U$. In particular, as $(z^1, \ldots, z^m) \rightarrow 0$,
$F_{\overline{x}}^* \omega_X$ approaches
$2 \sqrt{-1} \sum_{i=1}^m c_i \frac{dz^i \wedge
d\overline{z}^i}{|z^i|^2 \log^2(|z^i|^{-2})} + 
G_{\overline{x}}^* \omega_{D_I^{int}}$. 
Combining with the uniform
heat kernel asymptotics of $e^{-\epsilon \triangle^{(1,1)}}$ on  $U$, 
one sees that there is some $K > 0$ so that for any $k \ge K$,
there is $L^1$-convergence in (\ref{l1}) on $O_k$.
This is a contradiction. 

It follows that
\begin{equation}
\lim_{\epsilon \rightarrow 0} \int_{\overline{X}}
\left( e^{-\epsilon \triangle^{(1,1)}} \omega_X \right)^j \wedge \alpha
= \int_{X}
\omega_X^j \wedge \alpha.
\end{equation}
Thus $[\omega_X^j] = [\omega_X]^j$ in 
$\HH^{(j,j)}(\overline{X}; \R)$. \\
(3). The same argument as in part (1) shows that
$\Ric(\omega_X)$ extends to an element 
$[\Ric(\omega_X)] \in \HH^{(1,1)}(\overline{X}; \R)$.
From the asymptotics in Definition \ref{standard},
$\frac{h_{K_X}}{\prod_{i=1}^k |\sigma_i|^2_{L_i} \log^2(|\sigma_i|^{-2}_{L_i})}$
extends to a continuous Hermitian metric $h_{K_{\overline{X}}}$ on
$K_{\overline{X}}$.
Now $h_{K_{\overline{X}}} \otimes \bigotimes_{i=1}^k h_{L_i}$ is
a Hermitian metric on $K_{\overline{X}} \otimes L_D$ and on $X$, there
is an equality of currents:
\begin{equation}
\Ric(\omega_X) = - \sqrt{-1} F(h_{K_X}) = 
- \sqrt{-1} F(h_{K_{\overline{X}}} \otimes \bigotimes_{i=1}^k h_{L_i})
+ \sqrt{-1} \partial \overline{\partial} \sum_{i=1}^k
\log \log^2(|\sigma_i|^{-2}_{L_i}).
\end{equation}
The extension by zero of $\sqrt{-1} \partial \overline{\partial} \sum_{i=1}^k
\log \log^2(|\sigma_i|^{-2}_{L_i})$ to $\overline{X}$ is a 
closed
$(1,1)$-current 
which is the image under $\sqrt{-1} \partial \overline{\partial}$
of $\sum_{i=1}^k
\log \log^2(|\sigma_i|^{-2}_{L_i}) \in L^1(\overline{X})$. Hence
$[\sqrt{-1} \partial \overline{\partial} \sum_{i=1}^k
\log \log^2(|\sigma_i|^{-2}_{L_i})]$ vanishes in $\HH^{(1,1)}(X; \R)$ and so
$[\Ric(\omega_X)] = - 2\pi [K_{\overline{X}} + D]$.\\
(4). Let $\omega_{\overline{X}}$ be an arbitrary smooth
K\"ahler form on $\overline{X}$. From \cite[Theorem 4.2]{demailly-paun}, it suffices
to show that $\int_Y [\omega_{\overline{X}}]^{p-j} 
[\omega_X]^{j} > 0$ for all irreducible analytic sets $Y$ and
all $0 \le j \le p$, where $\dim_\C Y = p$.

For any $\epsilon > 0$,
$\int_Y [\omega_{\overline{X}}]^{p-j} 
[\omega_X]^{j} = \int_Y \omega_{\overline{X}}^{p-j} \wedge
\left( e^{-\epsilon \triangle^{(j,j)}} \omega_X^{j} \right)$.
Suppose first that $Y \not\subset D$. Then $Y-D$ is dense in $Y$.
By taking $\epsilon$ small
and using the asymptotics in Definition \ref{standard}, it follows that
$\int_Y [\omega_{\overline{X}}]^{p-j} [\omega_X]^{j} =
\int_Y \omega_{\overline{X}}^{p-j} \wedge \omega_X^{j} > 0$.
Now suppose that
$Y \subset D_I$ and $D_I$ is minimal with respect to this property.
For $x \in D_I^{int}$, let $R_x^* : \Lambda^{(j,j)}_x \overline{X} 
\rightarrow \Lambda_x^{(j,j)} D_I^{int}$ be the pullback map.
Using the asymptotics in Definition \ref{standard}, 
if $x \in D_I^{int}$ then
\begin{align}
\lim_{\epsilon \rightarrow 0} 
R_x^* \left( \left( e^{-\epsilon \triangle^{(j,j)}} \omega_X^{j} \right)(x)
\right) & =
\lim_{\epsilon \rightarrow 0} 
\int_{\overline{X}} 
R_x^* \left( e^{-\epsilon \triangle^{(j,j)}}(x,y) \: \omega_X^{j}(y) \right)
\: \dvol_{\overline{X}}(y) \\
& = \: 
\lim_{\epsilon \rightarrow 0} 
\int_{\overline{X}} 
(4 \pi \epsilon)^{-n} \: e^{- 
\frac{d_{\overline{X}}(x,y)^2}{4\epsilon}} R_x^* P_{x,y}
\omega_X^{j}(y) \:
\dvol_{\overline{X}}(y) \: = \: 
\omega_{D_I^{int}}^{j}(x). \notag
\end{align}
It follows that
$\int_Y [\omega_{\overline{X}}]^{p-j} [\omega_X]^{j} =
\int_Y \omega_{\overline{X}}^{p-j} \wedge \omega_{D_I^{int}}^{j} > 0$.
\end{proof}

\begin{remark}
Part (2) of Theorem \ref{extension} has a more direct proof if
$\omega_X$ has superstandard spatial asymptotics in the sense of
Definition \ref{superstandard} below.
Part (3) of Theorem \ref{extension} also follows from
\cite[\S 1]{Mumford}. 
\end{remark}

\begin{example} \label{stanexample}
Let ${\omega}_{\overline{X}}$ be a K\"ahler metric on
$\overline{X}$. Given positive numbers $\{c_i\}_{i=1}^k$, define a
$(1,1)$-form on $X$ by
\begin{align} \label{mess}
\omega_X & = {\omega}_{\overline{X}} - \sqrt{-1} \partial \overline{\partial}
\sum_{i=1}^k c_i \log \log^2 |\sigma_i|^{-2}_{L_i} \\
& =
{\omega}_{\overline{X}} - 2 \sqrt{-1} 
\frac{\partial \overline{\partial} \sum_{i=1}^k c_i
\log |\sigma_i|^{-2}_{L_i}}{\log |\sigma_i|^{-2}_{L_i}}
+ 2 \sqrt{-1} \sum_{i=1}^k 
c_i \frac{\partial \log |\sigma_i|^{-2}_{L_i}}{\log |\sigma_i|^{-2}_{L_i}} 
\wedge 
\frac{\overline{\partial} \log |\sigma_i|^{-2}_{L_i}}{\log |\sigma_i|^{-2}_{L_i}}. \notag
\end{align}
Now $\widehat{\omega} = {\omega}_{\overline{X}} 
+ 2 \sqrt{-1} \sum_{i=1}^k 
c_i \frac{\partial \log |\sigma_i|^{-2}_{L_i}}{\log |\sigma_i|^{-2}_{L_i}} 
\wedge 
\frac{\overline{\partial} \log |\sigma_i|^{-2}_{L_i}}{\log |\sigma_i|^{-2}_{L_i}}$ is
a K\"ahler metric on $X$.
For any $\epsilon > 0$, if the Hermitian metrics $h_{L_i}$ are multiplied by
a sufficiently small constant then
\begin{equation}
- \epsilon \widehat{\omega} \le  2 \sqrt{-1} 
\frac{\partial \overline{\partial} \sum_{i=1}^k c_i
\log |\sigma_i|^{-2}_{L_i}}{\log |\sigma_i|^{-2}_{L_i}} \le \epsilon 
\widehat{\omega}.
\end{equation}
Hence by rescaling the Hermitian metrics, we can achieve that
$\omega_X$ defines a K\"ahler metric on $X$. 

One can check that
$\omega_X$ has standard spatial asymptotics.
To describe $\omega_{D_I^{int}}$, suppose that
$D_I = \bigcap_{j=1}^m D_{i_j}$. After permuting indices, we
can assume that $D_I = D_1 \cap D_2 \cap \ldots \cap D_m$. Then
\begin{equation}
\omega_{D_I^{int}} = {\omega}_{\overline{X}} \Big|_{D_I^{int}} 
- \sqrt{-1} \partial \overline{\partial}
\sum_{i=m+1}^k c_i \log \log^2 |\sigma_i|^{-2}_{L_i},
\end{equation}
where the last computation is performed on $D_I^{int}$.
\end{example}

\begin{example} \label{ke}
Suppose that $\overline{X}$ is a compact K\"ahler manifold,
$D$ is an effective divisor in $\overline{X}$ with simple
normal crossings and $[K_{\overline{X}} + D] > 0$. We
use the K\"ahler-Ricci flow and
Theorem \ref{th:non-deg-convergence1}, or
\cite[Theorem 1.1]{chau}, 
to construct the K\"ahler-Einstein
metric on $X = \overline{X} - D$ which is known to exist from
\cite{kobayashi,tian-yau,tsu}.

The first step, as in \cite{kobayashi,tian-yau,tsu}, is to
construct a $0$-th order approximation to the K\"ahler-Einstein
metric using an idea of Carlson-Griffiths 
\cite[Proposition 2.1]{Carlson-Griffiths}.
Namely, since $[K_{\overline{X}} + D] > 0$, we can find a Hermitian
metric $h_{K_{\overline{X}} \otimes L_D}$ on
$K_{\overline{X}} \otimes L_D$ so that
$\sqrt{-1} F(h_{K_{\overline{X}} \otimes L_D}) > 0$. Fix
$\omega_{\overline{X}} = 
\sqrt{-1} F(h_{K_{\overline{X}} \otimes L_D})$.
Now perform the construction of Example \ref{stanexample} with
$c_1 = c_2 = \ldots = c_k = 1$ to get a K\"ahler metric $\omega_X$
on $X$, with corresponding Hermitian metric $h_{K_X}$ on $K_X$.
The construction also produces a Hermitian metric  $h_D$ on $L_D$.
This,
along with $h_{K_{\overline{X}} \otimes L_D}$, gives a Hermitian metric
$h_{K_{\overline{X}}}$ on $K_{\overline{X}}$.
Then 
\begin{equation}
- \Ric(\omega_X) + \sqrt{-1} \partial \overline{\partial}
\log \frac{h_{K_X}}{h_{K_{\overline{X}}}
\prod_{i=1}^k |\sigma_i|_{L_i}^2 \log^2 |\sigma_i|_{L_i}^{-2}}
\: = \: \omega_X
\end{equation}
on $X$.
From the standard spatial asymptotics,
$\log \frac{h_{K_X}}{h_{K_{\overline{X}}}
\prod_{i=1}^k |\sigma_i|_{L_i}^2 \log^2 |\sigma_i|_{L_i}^{-2}}$
has bounded covariant derivatives 
(with respect to $\omega_X$).
We can now apply Theorem \ref{th:non-deg-convergence1} or
\cite[Theorem 1.1]{chau} to $\omega_X$. 

However, to be more general,
suppose that $f_1$ is any smooth function on $X$ so that
\begin{enumerate}
\item $f_1$ has bounded covariant derivatives (with respect to
$\omega_X$) and
\item $\omega_X + \sqrt{-1} \partial \overline{\partial} f_1$
is a K\"ahler metric which is biLipschitz equivalent to $\omega_X$.
\end{enumerate}
Then
\begin{align}
& - \Ric(\omega_X + \sqrt{-1} \partial \overline{\partial} f_1) 
+ \\
& \sqrt{-1} \partial \overline{\partial}
\left(
\log \frac{\omega_X^n h_{K_X}}{
(\omega_X + \sqrt{-1} \partial \overline{\partial} f_1)^n
h_{K_{\overline{X}}}
\prod_{i=1}^k |\sigma_i|_{L_i}^2 \log^2 |\sigma_i|_{L_i}^{-2}} + f_1 \right)
\: = \notag \\
& \omega_X + \sqrt{-1} \partial \overline{\partial} f_1. \notag
\end{align}
Putting
\begin{equation}
f = \log \frac{\omega_X^n h_{K_X}}{
(\omega_X + \sqrt{-1} \partial \overline{\partial} f_1)^n
h_{K_{\overline{X}}}
\prod_{i=1}^k |\sigma_i|_{L_i}^2 \log^2 |\sigma_i|_{L_i}^{-2}} + f_1,
\end{equation}
Theorem \ref{th:non-deg-convergence1}, or
\cite[Theorem 1.1]{chau}, 
implies that the normalized K\"ahler-Ricci flow starting with the initial
metric
$\omega_0 = \omega_X  + \sqrt{-1} \partial \overline{\partial} f_1$
converges to a complete K\"ahler-Einstein
metric on $X$ with Einstein constant $-1$. 
(Such a metric is necessarily unique.)
From the evolution formulas for
the volume and scalar curvature under Ricci flow, one easily shows
that the K\"ahler-Einstein metric has finite volume.
In the case of complex dimension one, we recover some of the
results of \cite{Ji-Mazzeo-Sesum}.
\end{example}

\section{Preservation of standard spatial asymptotics}

In this section we show that the property of having standard spatial
asymptotics is preserved by the K\"ahler-Ricci flow.  We use this to
give an upper bound on the singularity time $T_{\sing}$.

\begin{theorem} \label{0thorder}
Suppose that $\omega_X(0)$ has standard spatial asymptotics associated to
$\{\omega_{D_I^{int}}(0)\}$ and
$\{c_i\}_{i=1}^k$. Suppose that the normalized K\"ahler-Ricci flow 
$\omega_X(t)$, with initial K\"ahler form $\omega_X(0)$,
exists on a maximal
time interval $[0, T)$ in the sense of Theorem \ref{th:existence1}. 
Then for all $t \in [0, T)$, $\omega_X(t)$ has standard
asymptotics associated to $\{\omega_{D_I^{int}}(t)\}$ and
$\{1 + (c_i -1) e^{-t} \}_{i=1}^k$, where $\omega_{D_I^{int}}(t)$ is a
normalized K\"ahler-Ricci flow solution on $D_I^{int}$ with initial
K\"ahler-form $\omega_{D_I^{int}}(0)$.
\end{theorem}
\begin{proof}
Suppose first that $D$ is a smooth divisor $C$ with a trivial
holomorphic normal bundle. Then there is a biholomorphic map
$F : \Delta \times C \rightarrow V$ to a neighborhood $V$ of $C$, with
$F$ restricting to the identity map from  ${\{0\} \times C}$ to $C \subset V$.
The restriction $F \Big|_{\Delta^* \times C} : \Delta^* \times C 
\rightarrow V \cap X$ has a
lift to a holomorphic covering map 
$\widetilde{F} : H \times C \rightarrow V \cap X$.  Suppose that the
conclusion of the theorem is not true.  Then for some $t^\prime \in [0, T)$
and some $\epsilon > 0$,
there is a sequence
$r_j \rightarrow \infty$ so that each
$\alpha_{r_j}^* \widetilde{F}^* \omega_X(t^\prime)$ has distance at least
$\epsilon$ from 
$(1 + (c_1 -1) e^{-t^\prime}) \frac{\sqrt{-1}}{2} \frac{dz^1 \wedge
d\overline{z}^1}{
(\text{Im}(z^1))^2}  + \omega_{C}(t^\prime)$
in the pointed $C^\infty$-topology.  (We use basepoint $\{\sqrt{-1}\} \times c_0$ for
some arbitrary $c_0 \in C$.) 

From our assumptions, there is a uniform positive lower bound on
the injectivity radius of $\widetilde{F}^* \omega_X(0)$ at
$\alpha_{r_j}(\{\sqrt{-1}\} \times c_0)$ or, equivalently, of
$\alpha_{r_j}^* \widetilde{F}^* \omega_X(0)$ at
$\{\sqrt{-1}\} \times c_0$.
By Hamilton's compactness theorem
\cite{Hamilton2},
after passing to a subsequence,
there is a pointed limit 
\begin{equation}
\lim_{j \rightarrow \infty}
\left( H \times C, \{\sqrt{-1}\} \times c_0,
\alpha_{r_j}^* \widetilde{F}^* \omega_X(t) \right) = 
\left( H \times C, \{\sqrt{-1}\} \times c_0,
\omega_\infty(t) \right)
\end{equation}
for some normalized
Ricci flow solution $\omega_\infty(t)$ on $H \times C$ that exists for
$t \in [0, T)$, with bounded curvature on compact time intervals. 
(Note that in taking the limit we do not have to perform diffeomorphisms.
Note also that the metric $\alpha_{r_j}^* \widetilde{F}^* \omega_X(t)$
on $H \times C$ is not complete, but nevertheless we can apply 
Ricci flow compactness to get the complete limiting metric
$\omega_\infty(t)$.)
Also by assumption,
$\lim_{j \rightarrow \infty}
\alpha_{r_j}^* \widetilde{F}^* \omega_X(0) =
c_1 \frac{\sqrt{-1}}{2} \frac{dz^1 \wedge d\overline{z}^1}{(\text{Im}(z^1))^2} + \omega_{C}(0)$. 
From the uniqueness of Ricci flow solutions with bounded curvature on
compact time intervals \cite{chen-zhu}, it follows that
$\omega_\infty(t^\prime) = 
(1 + (c_1 -1) e^{-t^\prime}) \frac{\sqrt{-1}}{2} \frac{dz^1 \wedge d\overline{z}^1}{
(\text{Im}(z^1))^2} + \omega_{C}(t^\prime)$, where
$\omega_C(\cdot)$ is a normalized Ricci flow solution on $C$.
This is a contradiction, thereby proving the theorem in this special
case. 

We now discuss the case when $D$ is a smooth divisor $C$ but its holomorphic
normal bundle need not be trivial.  In this case we may not be able to
use the covering space argument from before; for example, if $V$ is
a tubular neighborhood of $C$ then
$V-C$ may be simply connected and we cannot increase the injectivity radius
by passing to a cover.  
On the other hand, in some sense this problem is irrelevant,
since we can localize the argument and 
parametrize a neighborhood $U \subset \overline{X}$ of $\overline{x} \in C$
by a holomorphic map $F_{\overline{x}} : \Delta^n \rightarrow 
\overline{X}$ with 
$F_{\overline{x}}(0) = \overline{x}$ and
$F_{\overline{x}}(\{0 \} \times \Delta^{n-1}) \subset C$.
Then we can consider
the pullback metric $\widetilde{F}_{\overline{x}}^* \omega_X(t)$ on
the cover $H \times \Delta^{n-1}$ of $\Delta^* \times \Delta^{n-1}$ and
try to run the previous argument.  However, there is a new problem
because the
limiting metric on $H \times \Delta^{n-1}$ would not be complete, whereas 
the uniqueness results are for complete metrics.
Again, this
problem is somewhat irrelevant, since we should be able to patch 
together the local parametrizations
$F_{\overline{x}} : \Delta^n \rightarrow 
\overline{X}$ as $\overline{x}$ varies over $C$
and thereby effectively pass
to the setting of complete metrics.
To do so, it is convenient
to use the language of \'etale groupoids. 
We use the notion of a Ricci flow on an \'etale groupoid, as explained
in \cite[Section 5]{Lott1} and \cite[Section 3]{Lott2}.

Let us first reformulate the earlier setting, when the holomorphic normal
bundle is trivial, in terms of \'etale groupoids.
Let $T_v$ denote translation in
$H$ by $v \in \R$. Then $\alpha_r^{-1} T_v \alpha_r =
T_{r^{-1} v}$.
It follows that $\alpha_r^* \widetilde{F}^* \omega_X(0)$
is invariant under translation by $2 \pi r^{-1} \Z$. Then
$\lim_{r \rightarrow \infty}  
\alpha_r^* \widetilde{F}_{\overline{x}}^* \omega_X(0)
\: = \: \frac{\sqrt{-1}}{2} \frac{dz^1 \wedge d\overline{z}^1}{
(\text{Im}(z^1))^2}  + 
\omega_C(0)$, where the right-hand side is invariant under translation
by $\lim_{r \rightarrow \infty} 2 \pi r^{-1} \Z = \R_\delta$ on $H$. 
Here $\R_\delta$ denotes the group $\R$ with the discrete topology.
Equivalently, the pointed limit 
\begin{equation}
\lim_{r \rightarrow \infty} \left( X, 
\widetilde{F}(\alpha_r(\sqrt{-1}), c_0),
\omega_X(0) \right) \cong
\lim_{r \rightarrow \infty} \left( \Delta^* \times C, 
(\pi(\alpha_r(\sqrt{-1})), c_0), F^* \omega_X(0) \right)
\end{equation}
exists
as a pointed Riemannian groupoid, whose underlying \'etale groupoid is
the cross-product groupoid
$(H \times C) \rtimes \R_\delta$, with the K\"ahler form on the space of
units $H \times C$ being $c_1 \frac{\sqrt{-1}}{2} \frac{dz^1 \wedge d\overline{z}^1}{
(\text{Im}(z^1))^2}  + 
\omega_C(0)$. Then the normalized 
K\"ahler-Ricci flow solution on the \'etale groupoid 
is given by the $\R_\delta$-invariant 
normalized K\"ahler-Ricci flow solution
$((1 + (c_1 -1) e^{-t})) \frac{\sqrt{-1}}{2} \frac{dz^1 \wedge d\overline{z}^1}{
(\text{Im}(z^1))^2}
+ \omega_{C}(t)$ on the space
of units $H \times C$.

In the case when the holomorphic normal bundle of $C$ is not trivial,
we take $\overline{x} \in C$ and choose a local parametrization
$F_{\overline{x}} : \Delta^n \rightarrow \overline{X}$ with
$F_{\overline{x}}(0) = \overline{x}$ and 
$F_{\overline{x}}(\{0\} \times \Delta^{n-1}) \subset C$. Then
the pointed limit
$\lim_{r \rightarrow \infty} \left( X, 
\widetilde{F}_{\overline{x}}(\alpha_r(\sqrt{-1}), 0),
\omega_X(0) \right)$
exists
as a pointed Riemannian groupoid, whose underlying \'etale groupoid is
the cross-product groupoid
$(H \times C) \rtimes \R_\delta$, with the K\"ahler form on the space of
units $H \times C$ being $c_1 \frac{\sqrt{-1}}{2} \frac{dz^1 \wedge d\overline{z}^1}{
(\text{Im}(z^1))^2}  + 
\omega_C(0)$. Again, the  normalized
K\"ahler-Ricci flow solution on the \'etale groupoid 
is given by the $\R_\delta$-invariant normalized K\"ahler-Ricci flow solution
$((1 + (c_1 -1) e^{-t})) \frac{\sqrt{-1}}{2} \frac{dz^1 \wedge d\overline{z}^1}{
(\text{Im}(z^1))^2}
+ \omega_{C}(t)$ on the space
of units $H \times C$, where $\omega_C(\cdot)$ is a normalized 
K\"ahler-Ricci flow solution on
$C$. 

Now the uniqueness argument of \cite{chen-zhu} extends to Ricci flow
solutions on \'etale groupoids.
Along with the compactness result for Ricci flow
solutions on \'etale groupoids
\cite[Theorem 1.4]{Lott1}, we can prove the theorem using a
contradiction argument as before.

Finally, in the case of general $D$,
suppose that $\overline{x} \in D_I^{int}$. Let
$F_{\overline{x}} : \Delta^n \rightarrow U$ be the holomorphic
parametrization near $\overline{x}$. 
Then the pointed limit
$\lim_{r \rightarrow \infty} \left( X, 
\widetilde{F}_{\overline{x}}(\alpha_r(\sqrt{-1},\ldots,\sqrt{-1}), 0),
\omega_X(0) \right)$
exists
as a pointed Riemannian groupoid, whose underlying \'etale groupoid is
the cross-product groupoid
$(H^m \times D_I^{int}) \rtimes \R_\delta^m$, 
with the K\"ahler form on the space of
units $H^m \times D_I^{int}$ being 
$\sum_{i=1}^m c_i \frac{\sqrt{-1}}{2} \frac{dz^i \wedge 
d\overline{z}^i}{
(\text{Im}(z^i))^2}  + 
\omega_{D_I^{int}}(0)$. 
The normalized K\"ahler-Ricci flow solution on the \'etale groupoid 
is given by the $\R_\delta^m$-invariant normalized K\"ahler-Ricci flow solution
$\sum_{i=1}^m ((1 + (c_i -1) e^{-t})) \frac{\sqrt{-1}}{2} \frac{dz^i \wedge d\overline{z}^i}{
(\text{Im}(z^i))^2} + \omega_{D_I^{int}}(t)$ on the space
of units $H^m \times D_I^{int}$, where
$\omega_{D_I^{int}}(t)$ is a complete normalized K\"ahler-Ricci 
flow solution on $D_I^{int}$.
The theorem now follows from a
contradiction argument as before. 
\end{proof}

\begin{remark}
It follows that under the hypotheses of Theorem \ref{0thorder}, the
normalized K\"ahler-Ricci flow 
exists on each $D_I^{int}$ for a time interval of at least
$[0, T)$, with bounded curvature on compact subintervals of $[0, T)$.
Note in this regard 
that Theorem \ref{maintheorem} is consistent with passing to
the divisor, in the sense that $\left( K_{\overline{X}} + L_{D_I} \right)
\big|_{D_I} =
K_{D_I}$, and if $c$ is a K\"ahler class on $\overline{X}$ then its
pullback to $D_I$ is a K\"ahler class on $D_I$.
\end{remark}

\begin{remark}
Continuing with the previous remark,
the divisor $D$ is itself a complex space in the sense of
\cite[Chapter II.5]{demailly}. This suggests that one should be able
to extend the results of this paper from the setting of pairs
$(\overline{X}, D)$ to the setting of
complex spaces $Y$, or some class thereof. 
For example, a standard K\"ahler metric on such a
complex space would consist of complete K\"ahler metrics on the strata
$Y_j - Y_{j-1}$ of $Y$ having ``standard'' spatial asymptotics as
one approaches (in $Y_j$) a substratum $Y_k$ of $Y_j$. 
\end{remark}

\begin{corollary} \label{uasymp}
Let $u_{D_I^{int}}(t) \in C^\infty(D_I^{int})$ be the time-$t$ solution of 
(\ref{eq:skrf2}) on $D_I^{int}$. Then the time-$t$ solution $u_X(t) \in C^\infty(X)$ of
(\ref{eq:skrf2}) on $X$ has standard spatial asymptotics associated to
$\{ \const_I(t) + u_{D_I^{int}}(t) \}$, where $\const_I(t)$ is spatially constant and only depends on the time $t$.
\end{corollary}
\begin{proof}
This follows from (\ref{ueqn1}) and Theorem \ref{0thorder}.
\end{proof}

\begin{corollary}
Suppose that $\omega_X(0)$ has
standard spatial asymptotics associated to $\{\omega_{D_I^{int}}(0) \}$ and
$\{c_i\}_{i=1}^k$. Let $T_1$ and $T_2$ be the same as in Theorem \ref{th:existence1}.

Let $T_3$ be the supremum (possibly infinite) of the numbers 
$T^\prime$ so that there is a smooth solution for $u$ in 
(\ref{eq:skrf2}) on the time interval $[0,T^\prime]$ such that
$\omega_t + \sqrt{-1} 
\partial \bar{\partial} u$ is a K\"ahler metric with 
standard spatial asymptotics associated to $\{\omega_{D_I^{int}}(t) \}$ and
$\{(1 + (c_i -1) e^{-t})\}_{i=1}^k$.

Let $T_4$ be the supremum (possibly infinite) of the numbers 
$T$ for which there is a function $F_{T} \in C^\infty(X)$, with
standard spatial asymptotics associated to $\{u_{D_I^{int}}(T) \}$, such
that $\omega_{T} + \sqrt{-1} \partial \bar{\partial} F_{T}$ is a
K\"ahler metric .

Then $T_1 = T_2 = T_3 = T_4$.
\end{corollary}
\begin{proof}
This follows from Theorem \ref{th:existence1} and Corollary \ref{uasymp}.
\end{proof}

\begin{corollary}
Under the hypotheses of Theorem \ref{0thorder},
the maximal existence time $T_{\sing} \in (0, \infty]$ of the normalized 
K\"ahler-Ricci flow on $X$ is bounded 
above by 
\begin{equation}
\sup \{t \in \R^+ : e^{-t} [\omega_X(0)] + 2\pi (1 - e^{-t}) 
[K_{\overline{X}} + D] \text{ lies in the K\"ahler cone of } \overline{X} \}.
\end{equation} 
\end{corollary}
\begin{proof}
Suppose that $T^\prime < T_{\sing}$. From Theorem \ref{extension} and the
normalized Ricci flow equation, 
\begin{equation}
\frac{d}{dt} [\omega_X(t)] = 2\pi [K_{\overline{X}} + D] - [\omega_X(t)]
\end{equation}
in $\HH^{(1,1)}(\overline{X}; \R)$.
Thus 
\begin{equation}
[\omega_X(T^\prime)] = e^{-T^\prime} [\omega_X(0)] + 2\pi (1 - e^{-T^\prime}) 
[K_{\overline{X}} + D].
\end{equation}
Also from Theorem \ref{extension}, $[\omega_X(T^\prime)]$ is a
K\"ahler class on $\overline{X}$. The corollary follows. 
\end{proof}

\section{Superstandard spatial asymptotics} \label{super}

In this section we introduce the notion of superstandard spatial
asymptotics.  We show that having this property is preserved under
the Ricci flow.  We then prove the first part of Theorem \ref{maintheorem}.

We first prove a lemma regarding the singular support of
the $\partial \overline{\partial}$-operator applied to certain functions.
In general, if $X = \overline{X} - D$ and $J \in C^\infty(X)$, let
$\overline{J}$ be a measurable extension of $J$ to $\overline{X}$.
Suppose that $\overline{J} \in L^1(\overline{X})$.
Note that since $D$ has measure zero,
this element of $L^1(\overline{X})$ is independent
of the particular measurable extension of $J$ to $\overline{X}$ that we choose.
For concreteness, we will use the extension by zero.

Consider the form
$\sqrt{-1}\partial \overline{\partial} J \in \Omega^{(1,1)}(X)$ on $X$.
If it has finite mass then
we can extend it by zero to $\overline{X}$, to obtain
the current $\sqrt{-1} \overline{\partial \overline{\partial} J}$.
Or we could consider the
$(1,1)$-current 
$\sqrt{-1} \partial \bar{\partial} \bar{J}$ on
$\overline{X}$. These two currents do not have to be the same.  For example,
if $\overline{X} = S^2$ and $D = \pt$, suppose that
$J \in C^\infty(S^2 - \pt)$ equals 
$\log |z|$ in a neighborhood of $\pt = \{0\}$. Then
$\sqrt{-1} \overline{\partial \bar{\partial} J}$ has no singular
support on $S^2$ and represents a nonzero class in $\HH^2(S^2; \R)$, whereas
the current $\sqrt{-1} \partial \bar{\partial} \bar J$ has
singular support at $\pt$ and vanishes in $\HH^2(S^2; \R)$.

The next lemma gives a sufficient condition on $J$ for the
two extensions to agree on $\overline{X}$.

\begin{lemma} \label{coho}
Let $\omega_{\overline{X}}$ be a smooth K\"ahler form on $\overline{X}$.
Given $J \in C^\infty(X)$,
suppose that $|J(x)| = o(\log \prod_{i=1}^k |\sigma_i|_{L_i}^{-2})$ 
as $x \rightarrow D$.
Suppose that $\sqrt{-1} \partial \overline{\partial} J$ has locally
finite mass on $\overline{X}$
and there is some $C > 0$ so
that when restricted to $X \subset \overline{X}$,
\begin{equation} \label{pluri} 
\sqrt{-1} \partial \overline{\partial} J \ge - C \omega_{\overline{X}}.
\end{equation} 
If
$\overline{J}$ denotes the extension of $J$
by zero to $\overline{X}$ then the current 
$\sqrt{-1} \partial \bar{\partial} \bar J$
on $\overline{X}$ equals
$\sqrt{-1} \overline{\partial \overline{\partial} J}$,
the extension by zero of the current 
$\sqrt{-1} \partial \overline{\partial} J$ on $X$.
\end{lemma}
\begin{proof}
For $\epsilon' \ge 0$, define 
$J_{\epsilon'} \in C^\infty(X)$ by
\begin{equation*}
J_{\epsilon'}(x) = J(x) - \epsilon'
\log \prod_{i =1}^k |\sigma_i(x)|_{L_i}^{-2}.
\: \: \: \: \: \: \: \: \: \: \: \: \: \: \: \:
\: \: \: \: \: \: \: \: \: \: \: \: \: \: \: \:
\: \: \: \: \: \: \: \: \: \: \: \: \: \: \: \: (8.3)
\end{equation*}
If $\epsilon' > 0$ then $J_{\epsilon'}$ is bounded above
on $X$. We can find 
a neighborhood 
$U_{\epsilon'}$ of $D$ so that $J_{\epsilon'}$ is 
almost-plurisubharmonic on $U_{\epsilon'} - D$ in the sense of
(8.2).
As in [11, Theorem I.(5.24)], there is an
extension $\bar J_{\epsilon'}$ of $J_{\epsilon'}$
to $\overline{X}$ which is almost-plurisubharmonic on 
$U_{\epsilon'}$. Now
$\sqrt{-1} \partial \bar{\partial} \bar J_{\epsilon'}$
is a $(1,1)$-current on $\overline{X}$
which is smooth on $\overline{X}-D$.
Using [11, Theorem I.(5.8)], 
$\sqrt{-1} \partial \bar{\partial} \bar J_{\epsilon'}$ 
is measurable on $U_{\epsilon'}$.

Let 
$\sqrt{-1} 
\overline{\partial \overline{\partial} J_{\epsilon^\prime}}$
denote the extension by zero, to $\overline{X}$, of the current 
$\sqrt{-1} \partial \overline{\partial} J_{\epsilon^\prime}$ on $X$. 
As in the proof of part (1) of Theorem 6.6, 
$\sqrt{-1} \overline{\partial \overline{\partial} 
J_{\epsilon^\prime}}$
is a closed $(1,1)$-current on $\overline{X}$. Putting
\begin{equation}
T_{\epsilon^\prime} = \sqrt{-1} \partial \bar{\partial} \bar 
J_{\epsilon^\prime} -
\sqrt{-1} \overline{\partial \overline{\partial} 
J_{\epsilon^\prime}}
\end{equation}
gives a 
closed nonnegative measurable current which is supported on $D$.
Now
\begin{equation}
\sqrt{-1} \partial \bar{\partial} \bar 
J -
\sqrt{-1} \overline{\partial \overline{\partial}J} 
 = T_{\epsilon'} + \epsilon' 
\left( \sqrt{-1} \partial \bar{\partial} 
\log \prod_{i =1}^k |\sigma_i(x)|_{L_i}^{-2} \right) \Big|_D
\end{equation}
as currents on $\overline{X}$.
Since $\epsilon'$ was an arbitrary positive number, it follows that
$\sqrt{-1} \partial \bar{\partial} \bar 
J -
\sqrt{-1} \overline{\partial \overline{\partial}J}$ is a 
closed nonnegative measurable current which is supported on $D$.
Now
\begin{equation}
\sqrt{-1} \partial \bar{\partial} \bar 
J -
\sqrt{-1} \overline{\partial \overline{\partial}J} 
 = T_{\epsilon,\epsilon'} + \epsilon' 
\left( \sqrt{-1} \partial \bar{\partial} 
\log \prod_{i =1}^k |\sigma_i(x)|_{L_i}^{-2} \right) \Big|_D,
\end{equation}
as currents on $\overline{X}$.
Since $\epsilon'$ was an arbitrary positive number, it follows that
$\sqrt{-1} \partial \bar{\partial} \bar 
J -
\sqrt{-1} \overline{\partial \overline{\partial}J}$ is a 
closed nonnegative measurable current which is supported on $D$.

By \cite[Corollary III.(2.14)]{demailly},
$\sqrt{-1} \partial \bar{\partial} \bar 
J -
\sqrt{-1} \overline{\partial \overline{\partial}J} 
 = \sum_{i=1}^k c_i \delta_{D_i}$ for
some nonnegative constants $\{c_i\}_{i=1}^k$. However, it is easy
to show that if $c_i$ is nonzero then
$J$ has a logarithmic singularity near $D_i$; see 
the
Green-Riesz formula \cite[Proposition I.(4.22a)]{demailly}
and \cite[Example III.(6.9)]{demailly}. This 
contradicts the assumption on $J$.
\end{proof}

To motivate the definition of superstandard
spatial asymptotics, we first prove
a result about the Ricci curvature of a metric with standard 
spatial asymptotics.

\begin{lemma}
If $\omega_X$ has standard spatial asymptotics then we can write 
\begin{equation} \label{Riccisuper}
\Ric(\omega_X) = \eta^\prime_{\overline{X}} - 
\sqrt{-1} \partial \overline{\partial}
\left( - \sum_{i=1}^k \log \log^2 |\sigma_i|_{L_i}^{-2} + H^\prime \right)
\end{equation}
on $X$, 
where 
\begin{itemize}
\item $\eta^\prime_{\overline{X}}$ is a smooth closed $(1,1)$-form
on $\overline{X}$  with $[\eta^\prime_{\overline{X}}]=
- [K_{\overline{X}} + D]$,
and 
\item $H^\prime \in C^\infty(X) \cap L^\infty(X)$.
\end{itemize}
\end{lemma}
\begin{proof}
Choose a Hermitian
metric $h_{K_{\overline{X}} \otimes L_D}$ on
$K_{\overline{X}} \otimes L_D$.
Along with
Hermitian metrics $\{h_{L_i}\}_{i=1}^k$ on $\{L_i\}_{i=1}^k$, we obtain a 
Hermitian metric $h_{K_{\overline{X}}}$ on ${K_{\overline{X}}}$. Then
\begin{equation}
\Ric(\omega_X) = 
- \sqrt{-1} F(h_{K_{\overline{X}} \otimes L_D}) - \sqrt{-1} \partial
\overline{\partial}
\left( - \sum_{i=1}^k \log \log^2 |\sigma_i|_{L_i}^{-2} + 
\log \frac{h_{K_{\overline{X}}}
\prod_{i=1}^k |\sigma_i|_{L_i}^2 \log^2 |\sigma_i|_{L_i}^{-2}}{h_{K_X}}
\right)
\end{equation}
on $X$.
Put $\eta^\prime_{\overline{X}} =
- \sqrt{-1} F(h_{K_{\overline{X}} \otimes L_D})$ and
$H^\prime = \log \frac{h_{K_{\overline{X}}}
\prod_{i=1}^k |\sigma_i|_{L_i}^2 \log^2 |\sigma_i|_{L_i}^{-2}}{h_{K_X}}$.
By the standard spatial asymptotics, $H^\prime \in L^\infty(X)$.
The lemma follows.
\end{proof}

\begin{remark}
It follows from elliptic estimates that for each $k \geq 0$,
the function $H^\prime$ has uniform bounds
on its $k$-th covariant derivatives; see the end of the proof of
Theorem \ref{end1}.
\end{remark}

Recall that Definition \ref{standard} of standard asymptotics involves some
parameters $\{c_i\}_{i=1}^k$.

\begin{definition} \label{superstandard}
A K\"ahler metric 
$\omega_X$ on $X$ has superstandard spatial asymptotics if it has
standard spatial asymptotics and one can write
\begin{equation} \label{omegasuper}
\omega_X = \eta_{\overline{X}} - \sqrt{-1} \partial \overline{\partial}
\left( \sum_{i=1}^k c_i \log \log^2 |\sigma_i|_{L_i}^{-2} + H \right)
\end{equation}
where
\begin{itemize}
\item $\eta_{\overline{X}}$ is a smooth
closed $(1,1)$-form on $\overline{X}$,
\item $h_{L_i}$ is a Hermitian metric on the line bundle $L_i$ and
\item $H \in C^\infty(X) \cap L^\infty(X)$.
\end{itemize}
\end{definition}

\begin{example}
If $\overline{X} = S^2$ and $D = \pt$, suppose that
in terms of a local coordinate $z$ near $\pt$, the metric takes
the form $\omega_X = 
- \sqrt{-1} \partial \overline{\partial}
\left( \log \log^2 |z|^{-2} + \log \log \log^2 |z|^{-2} \right)$.
Then $\omega_X$ has standard asymptotics but does not have
superstandard asymptotics.
\end{example}

\begin{lemma} \label{indepofmetric}
The property of having superstandard spatial asymptotics is independent of the
choice of Hermitian metrics $\{h_{L_i} \}_{i=1}^k$.
\end{lemma}
\begin{proof}
Given a Hermitian metric $h_{L_i}$, any other Hermitian metric on $L_i$ can
be written as $\phi_i h_{L_i}$ for some positive $\phi_i \in C^\infty(\overline{X})$.Then
\begin{equation}
\log \log^2 (\phi_i^{-1} |\sigma_i|_{L_i}^{-2}) -
\log \log^2 |\sigma_i|_{L_i}^{-2} = 
2 \log \left( 1 + \frac{\log \phi_i^{-1}}{\log |\sigma_i|_{L_i}^{-2}}\right),
\end{equation}
which is bounded on $\overline{X}$. The lemma follows.
\end{proof}

\begin{example} \label{superstanexample}
Continuing with Example \ref{stanexample},
one can check that
$\omega$ has superstandard spatial asymptotics.
\end{example}

\begin{lemma} \label{same}
If $\omega_X$ has superstandard spatial asymptotics then 
$[\omega_X] = [\eta_{\overline{X}}]$ in 
$\HH^{(1,1)}(\overline{X}; \R)$.
\end{lemma}
\begin{proof}
Let $\omega_{\overline{X}}$ be a smooth K\"ahler form on $\overline{X}$.
First, 
from (\ref{omegasuper}) and the definition of standard asymptotics,
if $\overline{C} > 0$ is sufficiently large then there is some
$C > 0$ so that 
$\sqrt{-1} \partial \overline{\partial} 
(H - \overline{C} \sum_{i=1}^k c_i \log \log^2 |\sigma_i|_{L_i}^{-2}) \ge - C 
\omega_{\overline{X}}$ on $X$. Lemma \ref{coho} implies that
the extension of
$\sqrt{-1} \partial \overline{\partial} 
(H - \overline{C} \sum_{i=1}^k c_i \log \log^2 |\sigma_i|_{L_i}^{-2})$ 
by zero to $\overline{X}$ vanishes in
$\HH^{(1,1)}(\overline{X}; \R)$. It also follows from
Lemma \ref{coho} that the extension of
$- \sqrt{-1} \partial \overline{\partial} 
\sum_i c_i \log \log^2 |\sigma_i|_{L_i}^{-2}$ vanishes in
$\HH^{(1,1)}(\overline{X}; \R)$; see (\ref{mess}). Thus
$[\omega_X] = [\eta_{\overline{X}}] \in
\HH^{(1,1)}(\overline{X}; \R)$.
\end{proof}

\begin{theorem} \label{superextends}
Suppose that $\omega_X(0)$ has superstandard spatial asymptotics. 
Suppose that the normalized 
K\"ahler-Ricci flow $\omega_X(t)$, with initial K\"ahler
metric $\omega_X(0)$, exists on a maximal time interval $[0, T)$ in the
sense of Theorem \ref{th:existence1}. Then for all $t \in [0, T)$, $\omega_X(t)$ has
superstandard spatial asymptotics.
\end{theorem}
\begin{proof}
Recall the definition of $\omega_t$ from (\ref{omegat}).
Applying (\ref{Riccisuper}) and (\ref{omegasuper}) to $\omega_X(0)$,
we can write
\begin{align}
\omega_X(t) = & \omega_t + \sqrt{-1} \partial \overline{\partial} u(t) \\
= & 
- \eta^\prime_{\overline{X}} + e^{-t} (\eta_{\overline{X}} +
\eta^\prime_{\overline{X}}) - \notag \\
& \sqrt{-1} \partial \overline{\partial}
\left( \sum_i (1 + e^{-t} (c_i - 1)) \log \log^2 |\sigma_i|^{-2}_{L_i}
- H^\prime + e^{-t} (H + H^\prime)
- u(t) \right). \notag
\end{align}
From Corollary \ref{uasymp}, $u(t) \in L^\infty(X)$.
\end{proof}

\begin{theorem} \label{end1}
Suppose that $\omega_X(0)$ has superstandard spatial asymptotics.
Then the maximal existence time $T \in (0, \infty]$ of the K\"ahler-Ricci flow on $X$,
in the sense of Theorem \ref{th:existence1}, 
equals
\begin{equation}
\sup \{t \in \R^+ : e^{-t} [\omega_X(0)] + 2\pi (1 - e^{-t}) 
[K_{\overline{X}} + D] \text{ lies in the K\"ahler cone of } \overline{X} \}.
\end{equation} 
\end{theorem}
\begin{proof}
From Theorem \ref{th:existence1}, it suffices to show that if
$e^{-t} [\omega_X(0)] + 2\pi (1 - e^{-t}) 
[K_{\overline{X}} + D]$ lies in the K\"ahler cone of $\overline{X}$ then
there is a function $F_t \in C^\infty(X)$ such that
\begin{enumerate}
\item $\omega_t + \sqrt{-1} \partial \overline{\partial} F_t$ is a K\"ahler
metric which is biLipschitz equivalent to $\omega_X(0)$, and
\item For each $k$, the $k$-th covariant derivatives of $F_t$ (with respect
to the initial metric $\omega_X(0)$) are uniformly bounded.
\end{enumerate}

Suppose that ${\omega}_{\overline{X}}$ is a K\"ahler metric on
$\overline{X}$ whose
class in $\HH^{(1,1)}(\overline{X}; \R)$ equals $e^{-t} [\omega_X(0)] + 
2\pi (1 - e^{-t}) 
[K_{\overline{X}} + D]$. We construct a K\"ahler metric
$\omega_X$ on $X$ as in Example \ref{superstanexample}, using
the constants $\{ 1 + e^{-t} (c_i-1) \}_{i=1}^k$.
We now write
\begin{equation} \label{Feqn}
\omega_X = \omega_t + \sqrt{-1} \partial \overline{\partial} F
\end{equation}
and show that we can solve for $F$.
That is, we show that we can solve
\begin{equation} \label{tosolve}
\sqrt{-1} \partial \overline{\partial} F = 
\omega_X + {\rm Ric}(\omega_X(0))-e^{-t}\left(\omega_X(0)
+{\rm Ric}(\omega_X(0))\right).
\end{equation}

From Lemma \ref{indepofmetric}, for the purposes of the proof
we can assume that the Hermitian metrics $h_{L_i}$ are the
same in the construction of $\omega_X$ and in the
superstandard behavior of $\omega_X(0)$. Let $\eta_{\overline{X}}$ and
$\eta^\prime_{\overline{X}}$ be the $(1,1)$-forms on $\overline{X}$
involved in the superstandard behavior of $\omega_X(0)$.
From (\ref{mess}), (\ref{Riccisuper}) and (\ref{omegasuper}),
we can write
\begin{align} \label{comb1}
\omega_X + {\rm Ric}(\omega_X(0))-e^{-t}\left(\omega_X(0)
+{\rm Ric}(\omega_X(0))\right) = & {\omega}_{\overline{X}} + \eta^\prime_{\overline{X}}
- e^{-t} (\eta_{\overline{X}} + \eta^\prime_{\overline{X}}) \\
& - \sqrt{-1} \partial \overline{\partial} (H^\prime - e^{-t} (H + H^\prime)). \notag
\end{align}
From Lemma \ref{same} and our assumption on 
$\omega_{\overline{X}}$, we know that $\omega_{\overline{X}}$ and
$-\eta^\prime_{\overline{X}}
+ e^{-t} (\eta_{\overline{X}} + \eta^\prime_{\overline{X}})$ both represent
the same class in $\HH^{(1,1)}(\overline{X}; \R)$, namely
$e^{-t} [\omega_X(0)] + 2\pi (1 - e^{-t}) 
[K_{\overline{X}} + D]$. Thus 
\begin{equation} \label{comb2}
{\omega}_{\overline{X}} + \eta^\prime_{\overline{X}}
- e^{-t} (\eta_{\overline{X}} + \eta^\prime_{\overline{X}}) =
\sqrt{-1} \partial \overline{\partial} f
\end{equation} 
for some $f \in C^\infty(\overline{X})$.

From (\ref{comb1}) and (\ref{comb2}),
we can solve (\ref{Feqn}) for some
$F \in L^\infty(X)$. From (\ref{Feqn}), the Laplacian 
$\triangle_{\omega_X(0)} F = \Tr \left( \omega_X(0)^{-1} \sqrt{-1}
\partial \overline{\partial} F \right)$
has bounded $k$-th covariant derivatives (with respect to $\omega_X(0)$)
for each $k$. By elliptic regularity
(where near the divisor we work on the covering spaces 
$H^m \times \Delta^{n-m}$, 
which have bounded geometry), we conclude 
that $F$ also
has bounded $k$-th covariant derivatives
for each $k$. This proves the theorem. 
\end{proof}

This finishes the proof of the first part of Theorem
\ref{maintheorem}.  Theorem \ref{maintheorem} is stated for
the unnormalized K\"ahler-Ricci flow (\ref{unnorm}) instead of
the normalized K\"ahler-Ricci flow (\ref{eq:krf}), so one has to
make the translation between the two.

\section{Singularity type}

In this section we give sufficient conditions for the
K\"ahler-Ricci flow on a quasiprojective manifold to have
a type-II singularity.  We give examples in which this
happens.

\begin{theorem}
Suppose that $\omega_X(t)$ is a K\"ahler-Ricci flow solution 
on a quasiprojective manifold $X = \overline{X} - D$, $D \neq \emptyset$,
whose
initial metric $\omega_X(0)$ has superstandard spatial asymptotics.
Suppose that the maximal existence time $T_{\sing}$, in the sense of
Theorem \ref{th:existence1}, is finite. Suppose that there is a
number $C > 0$ 
so that for all $t \in [0, T_{\sing})$, we have $
\vol(X, g(t)) = \frac{1}{n!} \int_X \omega_X(t)^n \le
C (T_{\sing}-t)^n$. Then the Ricci flow has a type-II singularity at
time $T_{\sing}$, i.e. 
$\limsup_{t \rightarrow T_{\sing}} \left( (T_{\sing}-t) \sup_{x \in X}
|\Rm(x,t)| \right) = \infty$. 
\end{theorem}
\begin{proof}
If the theorem is not true then there is some $C^\prime > 0$ so that
for all $x \in X$ and $t \in [0, T_{\sing})$, we have
$|\Rm(x,t)| \le \frac{C^\prime}{T_{\sing}-t}$. From \cite[Theorem 1.4]{Naber},
for any $x^\prime \in X$ there is a sequence of times 
$t_i \rightarrow T_{\sing}$ so that
if we put $\tau_i = T_{\sing} - t_i$ then the rescaled Ricci flow solutions 
$g_i(x,t) = \tau_i^{-1} g \left( x, T_{\sing} + t \tau_i \right)$
have a pointed limit
$(X, g_i, (x^\prime,-1)) \stackrel{i \rightarrow \infty}{\rightarrow} 
(Y, g_\infty, (y_\infty, -1))$.
Here
$(Y, g_\infty)$ is a complete gradient shrinking soliton with bounded
curvature which is $\kappa$-noncollapsed at all scales, for some $\kappa >0$,
in the sense of Perelman \cite{Perelman}.
(Note that there is no $\kappa > 0$ so that the initial metric 
$\omega_X(0)$ is $\kappa$-noncollapsed at all scales. Nevertheless, 
in this setting the
blowup limit is $\kappa$-noncollapsed at all scales for some $\kappa$;
see \cite[Remark 2.2]{Naber}.)

From our assumption, $(Y, g(-1))$ has finite volume. However, the
$\kappa$-noncollapsing now implies that $Y$ is compact. 
(We thank Lei Ni for this remark.) Namely, if $Y$ is noncompact
then it contains an infinite sequence of disjoint unit balls.
The $\kappa$-noncollapsing, along with the bounded curvature, 
implies that there is a uniform positive lower bound on the
volumes of these balls.  This is a contradiction.

Thus $Y$ is compact. This implies that $X$ is compact, which
is a contradiction. The theorem follows.
\end{proof}

\begin{corollary}
Suppose that $\omega_X(t)$ is a K\"ahler-Ricci flow solution 
on a quasiprojective manifold $X = \overline{X} - D$, $D \neq \emptyset$,
whose
initial metric $\omega_X(0)$ has superstandard spatial asymptotics.
If $T_{\sing} < \infty$ and
$\lim_{t \rightarrow T_{\sing}} [\omega_X(t)] = 0$ in
$\HH^{(1,1)}(\overline{X}; \R)$ then there is a type-II singularity
at time $T_{\sing}$.
\end{corollary}
\begin{proof}
From the smoothness of $[\omega_X(t)]$, we can write
$[\omega_X(t)] = (T_{\sing} - t) R(t)$ for some smooth function
$R : [0, T_{\sing}] \rightarrow \HH^{(1,1)}(\overline{X}; \R)$.
Then there is a constant $C < \infty$ so that for 
$t \in [0, T_{\sing})$,
\begin{equation}
\int_X \omega_X^n(t) = \int_{\overline{X}} [\omega_X(t)]^n \le
C  (T_{\sing} - t)^n.
\end{equation}
The corollary follows.
\end{proof}

This finishes the proof of the second part of Theorem
\ref{maintheorem}.  We now give some examples, using the
unnormalized K\"ahler-Ricci flow of Theorem \ref{maintheorem}.

\begin{example} \label{surfaceexample}
Suppose that $\overline{X} = S^2$ and $D = \pt$, so
$X = S^2 - \pt = \R^2$. Let
$[S^2] \in \Image \left( \HH^2(S^2; \Z) \rightarrow \HH^2(S^2; \R)\right)
\cap \HH^{(1,1)}(S^2; \R)$
denote the fundamental class in cohomology.
Then $[K_{\overline{X}}] = - 2 [S^2]$ and $[D] = [S^2]$.
From Theorem \ref{maintheorem}, $T_{\sing}$ is
the supremum of
the numbers $T>0$ so that $[\omega_0] - 2\pi T [S^2] \in 
\HH^{(1,1)}(\overline{X}; \R)$ is a K\"ahler class. That is,
$T_{\sing} = \frac{1}{2\pi} \int_{\R^2} \omega_X(0) = 
\frac{1}{2\pi} \Vol(\R^2, g(0))$. (As we are now dealing with the unnormalized
K\"ahler-Ricci equation $\frac{d\omega}{dt} = - \Ric$, the
singularity time given here differs by a factor of two from the result
$\frac{1}{4\pi} \Vol(\R^2, g(0))$
stated in the introduction for the unnormalized Ricci flow
$\frac{dg}{dt} = - 2 \Ric$.)

As $[\omega_0] - 2\pi T_{\sing} [S^2]$ vanishes in  
$\HH^{(1,1)}(\overline{X}; \R)$,
we conclude that there is a type-II singularity at time $T_{\sing}$,
in agreement with the results of
Daskalopoulos-del Pino-Hamilton-Sesum 
\cite{Daskalopoulos-del PinoI,Daskalopoulos-del PinoII,Daskalopoulos-Hamilton,Daskalopoulos-Sesum}.
\end{example}

\begin{example}
Taking a product of the previous example with $S^2$, 
suppose that $\overline{X} = S^2 \times S^2$ and $D = \{\pt\} \times S^2$.
Let 
$[S^2]_1, [S^2]_2 \in \Image \left( \HH^2(S^2 \times S^2; \Z) \rightarrow
\HH^2(S^2 \times S^2; \R) \right) \cap 
\HH^{(1,1)}(S^2 \times S^2; \R)$
denote the fundamental classes of the two sphere factors.  Then
$K_{\overline{X}} = -2[S^2]_1 - 2 [S^2]_2$ and $D = [S^2]_2$.
We conclude that $T_{\sing}$ is the supremum of the times $T$ so that
$\int_{[S^2]_1} \omega_X(0) - 4 \pi T > 0$ and 
$\int_{[S^2]_2} \omega_X(0) - 2 \pi T > 0$. 
\begin{itemize}
\item If $\frac{\int_{[S^2]_1} \omega_X(0)}{4\pi} < 
\frac{\int_{[S^2]_2} \omega_X(0)}{2\pi}$ then $T_{\sing} =
\frac{\int_{[S^2]_1} \omega_X(0)}{4\pi}$. Since
$[\omega_X(0)] + 2\pi T_{\sing} [K_{\overline{X}}+D]$ is nonvanishing,
we cannot conclude that there is a type-II singularity.  In fact,
if the initial metric $\omega_X(0)$ is a product metric then the
first $S^2$-factor shrinks to a point at the singularity time before the
other factor can collapse, and we have a type-I singularity.
\item If $\frac{\int_{[S^2]_1} \omega_X(0)}{4\pi} = 
\frac{\int_{[S^2]_2} \omega_X(0)}{2\pi}$ then $T_{\sing}$ is this
common value. Since
$[\omega_X(0)] + 2\pi T_{\sing} [K_{\overline{X}}+D] = 0$, there is a
type-II singularity.
\item If $\frac{\int_{[S^2]_1} \omega_X(0)}{4\pi} > 
\frac{\int_{[S^2]_2} \omega_X(0)}{2\pi}$ then $T_{\sing} =
\frac{\int_{[S^2]_2} \omega_X(0)}{2\pi}$. Since
$[\omega_X(0)] + 2\pi T_{\sing} [K_{\overline{X}}+D]$ is nonvanishing,
we cannot conclude that there is a type-II singularity, although
there is one if $\omega_X(0)$ is a product metric.
\end{itemize}
\end{example}

\begin{example}
Suppose that $\overline{X} = \C P^n$ and $D$ consists of $k$ copies of
$\C P^{n-1}$ in general position. Let 
$[H] \in \Image \left( \HH^2(\C P^n; \Z) \rightarrow 
\HH^2(\C P^n; \R) \right) \cap
\HH^{(1,1)}(\C P^n; \R)$
be the hyperplane class.
Then $[K_{\overline{X}}] = - (n+1) [H]$ and $[D] = k [H]$, so
$T_{\sing}$ is
the supremum of
the numbers $T>0$ so that $[\omega_0] + 2\pi (-n-1+k) T [H] \in 
\HH^{(1,1)}(\overline{X}; \R)$ is a K\"ahler class.
\begin{itemize}
\item If $k > n+1$ then $T_{\sing} = \infty$. In this case there is
a finite-volume K\"ahler-Einstein metric $\omega_{KE}$ on $X$ with
Einstein constant $-1$
\cite{kobayashi,tian-yau,tsu}.
Theorem \ref{th:non-deg-convergence1} 
says that for a wide class of initial metrics,
the normalized K\"ahler-Ricci flow will converge to $\omega_{KE}$.
\item If $k=n+1$ then $T_{\sing} = \infty$. In this case there is
a complete Ricci-flat K\"ahler metric $\omega_{Ricci-flat}$ on $X$
\cite{tian-yau2}.
It should be possible to show
that for a large class of initial metrics, the unnormalized K\"ahler-Ricci
flow converges to a multiple of $\omega_{Ricci-flat}$.
\item If $k < n+1$ then $T_{\sing} = 
\frac{\int_{\C P^1 \cap X} \omega_X(0)}{2\pi(n+1-k)}$,
where $\C P^1$ denotes a generic complex line in $\overline{X} = \C P^n$.
If in addition $k \neq 0$ then there is a type-II singularity.

Note
that when $k=1$, there is a $U(n)$-invariant superstandard
initial
K\"ahler metric on $X = \C^n = \C P^n - \C P^{n-1}$. 
At infinity, it looks like a family of hyperbolic cusps parametrized by
$\C P^{n-1}$. It is plausible that in this case,
there is a rescaling limit at the singular time which is 
a $U(n)$-invariant gradient steady K\"ahler-Ricci soliton. Examples of
the latter are in \cite{Cao2}. From 
\cite{Daskalopoulos-del PinoI,Daskalopoulos-del PinoII,Daskalopoulos-Hamilton,Daskalopoulos-Sesum}, we know that there is such a rescaling limit when
$n=1$.
\end{itemize}
\end{example}

\end{document}